\documentclass[12p]{article}
\def\leaderdot{\leaders\hbox to 1 em {\hss.\hss}\hfill}

\def\skiplines#1 { \dimen3=\dimen2 \multiply\dimen3 by #1 \vskip \dimen3}
\def\fullline{\hbox to \fullhsize}

\def\numpage{\baselineskip=24pt\fullline{\the\footline}}

\def\mathcedilla{\vtop{\hbox{c}{\kern0pt\nointerlineskip}
             {\hbox{$\mkern-2mu \mathchar"0018\mkern-2mu$}}}}

\mathchardef\gq="0060 \mathchardef\dq="0027
\mathchardef\ssmath="19 \mathchardef\aemath="1A
\mathchardef\oemath="1B \mathchardef\omath="1C
\mathchardef\AEmath="1D \mathchardef\OEmath="1E
\mathchardef\Omath="1F \mathchardef\imath="10
\mathchardef\fmath="0166 \mathchardef\gmath="0167
\mathchardef\vmath="0176
\def\Irm{{\rm I}}
\def\IIrm{{\rm II}}
\def\IIIrm{{\rm III}}

\def\Vrm{{\rm V}}

\def\colleft{\strut\kern.3em}
\def\colright{\kern0pt}

\def\figureh{\hbox to}

\def\m@th{\mathsurround=0pt}
\newif\ifdtpt
\def\displ@y{\openup1\jot\m@th
    \everycr{\noalign{\ifdtpt\dt@pfalse
    \vskip-\lineskiplimit \vskip\normallineskiplimit
    \else \penalty\interdisplaylinepenalty \fi}}}
\def\eqalignl#1{\,\vcenter{\openup1\jot\m@th
                \ialign{\strut$\displaystyle{##}$\hfil&
                              $\displaystyle{{}##}$\hfil&
                              $\displaystyle{{}##}$\hfil&
                              $\displaystyle{{}##}$\hfil&
                              $\displaystyle{{}##}$\hfil\crcr#1\crcr}}\,}
\def\eqalignnol#1{\displ@y\tabskip\centering \halign to \displaywidth{
                  $\displaystyle{##}$\hfil\tabskip=0pt &
                  $\displaystyle{{}##}$\hfil\tabskip=0pt &
                  $\displaystyle{{}##}$\hfil\tabskip=0pt &
                  $\displaystyle{{}##}$\hfil\tabskip=0pt &
                  $\displaystyle{{}##}$\hfil\tabskip\centering &
                  \llap{$##$}\tabskip=0pt \crcr#1\crcr}}
\def\leqalignnol#1{\displ@y\tabskip\centering \halign to \displaywidth{
                   $\displaystyle{##}$\hfil\tabskip=0pt &
                   $\displaystyle{{}##}$\hfil\tabskip=0pt &
                   $\displaystyle{{}##}$\hfil\tabskip=0pt &
                   $\displaystyle{{}##}$\hfil\tabskip=0pt &
                   $\displaystyle{{}##}$\hfil\tabskip\centering &
                   \kern-\displaywidth\rlap{$##$}\tabskip=\displaywidth
                   \crcr#1\crcr}}
\def\eqalignc#1{\,\vcenter{\openup1\jot\m@th
                \ialign{\strut\hfil$\displaystyle{##}$\hfil&
                              \hfil$\displaystyle{{}##}$\hfil&
                              \hfil$\displaystyle{{}##}$\hfil&
                              \hfil$\displaystyle{{}##}$\hfil&
                              \hfil$\displaystyle{{}##}$\hfil\crcr#1\crcr}}\,}
\def\eqalignnoc#1{\displ@y\tabskip\centering \halign to \displaywidth{
                  \hfil$\displaystyle{##}$\hfil\tabskip=0pt &
                  \hfil$\displaystyle{{}##}$\hfil\tabskip=0pt &
                  \hfil$\displaystyle{{}##}$\hfil\tabskip=0pt &
                  \hfil$\displaystyle{{}##}$\hfil\tabskip=0pt &
                  \hfil$\displaystyle{{}##}$\hfil\tabskip\centering &
                  \llap{$##$}\tabskip=0pt \crcr#1\crcr}}
\def\leqalignnoc#1{\displ@y\tabskip\centering \halign to \displaywidth{
                  \hfil$\displaystyle{##}$\hfil\tabskip=0pt &
                  \hfil$\displaystyle{{}##}$\hfil\tabskip=0pt &
                  \hfil$\displaystyle{{}##}$\hfil\tabskip=0pt &
                  \hfil$\displaystyle{{}##}$\hfil\tabskip=0pt &
                  \hfil$\displaystyle{{}##}$\hfil\tabskip\centering &
                  \kern-\displaywidth\rlap{$##$}\tabskip=\displaywidth
                  \crcr#1\crcr}}
\def\doubleup#1{\,\vbox{\ialign{\hfil$##$\hfil\crcr
                  \mathstrut #1 \crcr }}\,}
\def\doublelow#1{\,\vtop{\ialign{\hfil$##$\hfil\crcr
                 \mathstrut #1 \crcr}}\,}
\def\charlvmidlw#1#2{\,\vtop{\ialign{##\crcr
      #1\crcr\noalign{\kern1pt\nointerlineskip}
      $\hfil#2\hfil$\crcr}}\,}
\def\charlvlowlw#1#2{\,\vtop{\ialign{##\crcr
      $\hfil#1\hfil$\crcr\noalign{\kern1pt\nointerlineskip}
      #2\crcr}}\,}
\def\charlvmidup#1#2{\,\vbox{\ialign{##\crcr
      $\hfil#1\hfil$\crcr\noalign{\kern1pt\nointerlineskip}
      #2\crcr}}\,}
\def\charlvupup#1#2{\,\vbox{\ialign{##\crcr
      #1\crcr\noalign{\kern1pt\nointerlineskip}
      $\hfil#2\hfil$\crcr}}\,}
\def\vspce{\kern4pt} \def\hspce{\kern4pt}    % rule spacing

\def\emptybox{\vbox{\kern.7ex\hbox{\kern.5em}\kern.7ex}}
\newfam\bfmitfam

\def\twodot{.\kern-0.1em.}

\def\paral{\mathrel{/\kern-.25em/}}
\def\grlo{\mathrel{\hbox{\lower.2ex\hbox{\rlap{$>$}\raise1ex\hbox{$<$}}}}}
\def\logr{\mathrel{\hbox{\lower.2ex\hbox{\rlap{$<$}\raise1ex\hbox{$>$}}}}}
\def\greq{\mathrel{\hbox{\lower1ex\hbox{\rlap{$=$}\raise1.2ex\hbox{$>$}}}}}
\def\loeq{\mathrel{\hbox{\lower1ex\hbox{\rlap{$=$}\raise1.2ex\hbox{$<$}}}}}
\def\grsim{\mathrel{\hbox{\lower1ex\hbox{\rlap{$\sim$}\raise1ex\hbox{$>$}}}}}
\def\losim{\mathrel{\hbox{\lower1ex\hbox{\rlap{$\sim$}\raise1ex\hbox{$<$}}}}}
\font\ninerm=cmr9
\def\uniset{\rlap{\ninerm 1}\kern.15em 1}

\def\emptysq{\mathbin{\vbox{\hrule\hbox{\vrule height1ex \kern.5em
                            \vrule height1ex}\hrule}}}
\def\emptyrect{\mathbin{\vbox{\hrule\hbox{\vrule height1ex \kern1em
                              \vrule height1ex}\hrule}}}
\def\rightonleftarrow{\mathrel{\hbox{\raise.5ex\hbox{$\rightarrow$}\ignorespaces
                                   \lower.5ex\hbox{\llap{$\leftarrow$}}}}}
\def\leftonrightarrow{\mathrel{\hbox{\raise.5ex\hbox{$\leftarrow$}\ignorespaces
                                   \lower.5ex\hbox{\llap{$\rightarrow$}}}}}
\def\bkB{{\rm I\kern-.17em B}}
\def\bkC{{\rm \kern.24em
            \vrule width.05em height1.4ex depth-.05ex
            \kern-.26em C}}
\def\bkD{{\rm I\kern-.17em D}}
\def\bkE{{\rm I\kern-.17em E}}
\def\bkF{{\rm I\kern-.17em F}}
\def\bkG{{\rm \kern.24em
            \vrule width.05em height1.4ex depth-.05ex
            \kern-.26em G}}
\def\bkH{{\rm I\kern-.22em H}}
\def\bkI{{\rm I\kern-.22em I}}
\def\bkJ{{\rm \kern.19em
            \vrule width.02em height1.5ex depth0ex
            \kern-.20em J}}
\def\bkK{{\rm I\kern-.22em K}}
\def\bkL{{\rm I\kern-.17em L}}
\def\bkM{{\rm I\kern-.22em M}}
\def\bkN{{\rm I\kern-.20em N}}
\def\bkO{{\rm \kern.24em
            \vrule width.05em height1.4ex depth-.05ex
            \kern-.26em O}}
\def\bkP{{\rm I\kern-.17em P}}
\def\bkQ{{\rm \kern.24em
            \vrule width.05em height1.4ex depth-.05ex
            \kern-.26em Q}}
\def\bkR{{\rm I\kern-.17em R}}
\def\bkT{{\rm \kern.24em
            \vrule width.02em height1.5ex depth 0ex
            \kern-.27em T}}
\def\bkU{{\rm \kern.30em
            \vrule width.02em height1.47ex depth-.05ex
            \kern-.32em U}}
\def\bkZ{{\rm Z\kern-.32em Z}}
\usepackage{graphicx}
\begin{document}
\centerline{{\bf SUR LES ACTIONS SYMPLECTIQUES QUASI-PRIMITIVES} }
\smallskip
\centerline{NGUIFFO B. BOYOM}
\smallskip
\centerline{UPRES A 5030}
\bigskip
{\bf Abstract}  A transitive smooth action of a connected Lie group $G$ on a manifold $M$ is called almost primitive (resp. primitive ) if $G$ doesn't contain any proper subgroup ( resp. any proper normal subgroup) whose induced action on $M$ is transitive as well. The aim of the present work is to investigate some combinatory properties of symplectic actions of completely solvable Lie groups. \par
\skiplines 2
{\bf Introduction}
Les vari\'et\'es en jeu ainsi que les objets consid\'er\'es sur ces vari\'et\'es sont lisses. Une loi d'op\'eration transitive d'un groupe de Lie connexe $G$ sur une vari\'et\'e $M$ est dite quasi-primitive (resp. primitive) si $G$ ne contient pas de sous-groupe propre (resp. sous-groupe distingu\'e propre) dont l'action induite sur $M$ est transitive. Ce travail est consacr\'e \`a  quelques aspects combinatoires des actions symplectiques des groupes de Lie connexes compl\`etement r\'esolubles. Si $G$ est un groupe de Lie connexe compl\`etement r\'esoluble on d\'esigne par ${\cal F}(G)$ l'ensemble des suites de composition de $G$ . Toute action symplectique transitive de $G$ d\'etermine une application de ${\cal F}(G)$ dans un ensemble de graphes ( resp. diagrammes ) analogues de graphes (resp. diagrammes ) de DYNKIN des groupes de Lie complexes semi-simples . Les caract\`eres quasi-primitif ou primitif de l'action symplectique correspondent \`a  des propri\`et\'es particuli\`eres de certains diagrammes. On peut y d\'eceler en particulier l'existence ou non de feuilletage symplectique invariant par $G$ ou celle de feuilletage Lagragien invariant par $G$ .\\
Ce travail a \'et\'e motiv\'e par la conjecture des tores plats dont l'\'enonc\'e est donn\'e ci-dessous. Soit $G$ un groupe de Lie connexe compl\`etement r\'esoluble contenant un r\'eseau cocompact $\Gamma$ .\\
\lq\lq Si la vari\'et\'e $ \Gamma \backslash G $ admet une m\'etrique Kahl\'erienne alors elle est diff\'eomorphe au tore plat.\rq\rq \\
Cette conjecture est vraie si $G$ est nilpotent. La premi\`ere d\'emonstration est due conjointement \`a C. BENSON et C. GORDON $\lbrack BG1 \rbrack$ . Une autre d\'emonstration a \'et\'e publi\'ee par Dusa McDUFF $\lbrack McD1 \rbrack$. Ces deux d\'emonstrations utilisent des propri\'et\'es des groupes de Lie nilpotents qui n'ont pas leur analogue dans les groupes de Lie r\'esolubles g\'en\'eraux. Les outils d\'ev\'elopp\'es dans ce travail ont pour but de contourner cette difficult\'e; ils sont utilis\'es dans $\lbrack NB6 \rbrack$ pour d\'emontrer une variante de la conjecture \'enonc\'ee ci-dessus et r\'epondre \`a d'autres questions pos\'ees dans $\lbrack BG2 \rbrack$ .\\
L'auteur remercie le r\'ef\'er\'e et un des \'editeurs de leur critiques et suggestions.
\skiplines 1
\noindent AMS classification. \par
\noindent Primary : 22E 25, 22E45, 53A10 \par
\noindent Secondary : 57S10 \par
\noindent Key words : bilagrangiens, diagrammes orient\'es, graphes,
r\'eductions Kahl\'eriennes. \par
\vfill\eject
\noindent{\bf 1.} {\bf GRAPHES ET DIAGRAMMES DES 2-FORMES FERM\'EES} \par
\skiplines 4
Les groupes de Lie consid\'er\'es sont r\'eels, connexes et simplement
connexes. Soit $ G $ un tel groupe, on entendra par alg\`ebre de Lie du groupe
de Lie $ G $ l'espace vectoriel $ g $ constitu\'e des champs de vecteurs
tangents invariants par les translations \`a gauche et muni du crochet de
Poisson des champs de vecteurs. On d\'esignera par $ \Omega^ 2_{\ell}( G) $
l'espace vectoriel
des 2-formes diff\'erentielles ferm\'ees qui sont invariantes par les
translations \`a gauche. \par
\bigskip
 {1.1. {\bf Graphes d'une 2-forme ferm\'ee.}} \\
Soit $ G $ un groupe de Lie compl\`etement r\'esoluble dont l'alg\`ebre de
Lie est not\'ee $ g\, . $ On d\'esigne par $ {\cal F}(G) $ l'ensemble des
suites de
composition dans $ G\, . $ Un \'el\'ement $ F \in  {\cal F}(G) $ est une
filtration de $ G $ par des
sous-groupes $ G_k\, , $ $ \leq  k \leq  \dim G\, , $ qui jouissent des
propri\'et\'es suivantes
$$ dim G_k = k; \leqno (i) $$
$$ G_k \hbox{ est un sous groupe distingu\hbox{\'e} dans } G_{k+1} \leqno (ii)$$
Une filtration $ F $ qui satisfait les propri\'et\'es (i) et (ii) ci-dessus sera dite {\it normale} lorsque chaque $ G_k $ est distingu\'e dans $ G\,; $ $ F $ est alors appel\'e suite de Jordan-H\"older de $ G\,; $ .
\noindent Un groupe compl\`etement r\'esoluble poss\`ede des suites de
compositions normales et
r\'eciproquement l'existence d'une suite de composition normale entra\^\i ne
la
compl\`ete r\'esolubilit\'e .\par
\smallskip
On consid\`ere $ \omega  \in  \Omega^ 2_{\ell}( G) $ et $ F \in  {\cal
F}(G)\, ; $ soit $ i_k: G_k\rightarrow  G $
l'homomorphisme inclusion de $ G_k $ dans $ G $ et soit $ \omega_ k= i^\ast_
k\omega \, . $ Si $ g_k $ est
l'alg\`ebre de Lie de $ G_k\, , $ on note $ h_k $ le sous-espace vectoriel des
$ \xi  \in  g_k $ solutions de $ i(\xi) \omega_ k= 0\, ; $ $ i(\xi) $ est le
produit int\'erieur par $ \xi \, . $ \par
\noindent Puisque $ \omega_ k $ est ferm\'ee, $ h_k $ est une sous-alg\`ebre
de Lie de l'alg\`ebre
de Lie $ g_k\, . $ On d\'esignera par $ H_k $ le sous-groupe de Lie connexe
associ\'e \`a la sous-alg\`ebre $ h_k\, . $ On associe ainsi au couple $
(\omega ,F) $ la
famille des triplets $ \left[G_k,H_k,\omega_ k \right]\, . $ Les sous-groupes
$ G_k\, , $ $ H_k $ sont ferm\'es dans
$ G\, , $ puisque ce dernier est simplement connexe. De la compl\`ete
r\'esolubilit\'e
de $ G\, , $ on d\'eduit l'existence de relation d'inclusion entre les deux
sous-groupes $ H_k $ et $ H_{k+1} $ (voir $\lbrack$NB 4$\rbrack$, $\lbrack$NB
5$\rbrack$). On a donc n\'ecessairement un
des deux diagrammes commutatifs d'homomorphismes de groupes de Lie
suivants
$$ \matrix{\displaystyle H_k & \displaystyle \rightarrow & \displaystyle
H_{k+1} \cr\displaystyle \downarrow & \displaystyle  & \displaystyle
\downarrow \cr\displaystyle G_k & \displaystyle \rightarrow & \displaystyle
G_{k+1}\, , \cr} \leqno (1) $$
$$ \matrix{\displaystyle H_k & \displaystyle \leftarrow & \displaystyle
H_{k+1} \cr\displaystyle \downarrow & \displaystyle  & \displaystyle
\downarrow \cr\displaystyle G_k & \displaystyle \rightarrow & \displaystyle
G_{k+1} \cr} \leqno (2) $$
Naturellement les fl\`eches sont des homomorphismes inclusion. \par
On d\'esigne par $ M_k $ le quotient $ G_k/H_k $ des classes \`a droite modulo
$ H_k $
des \'el\'ements de $ G_k\, ; $ $ M_k $ h\'erite de $ \left[G_k,H_k,\omega_ k
\right] $ d'une 2-forme symplectique
qu'on notera encore $ \omega_{ k\, .} $ Naturellement $ \left(M_k,\omega_ k
\right) $ est une vari\'et\'e
symplectique homog\`ene sous l'action \`a gauche de $ G_k\, . $ Lorsque $ H_k
$ est inclu
dans $ H_{k+1}\, , $ le sous-groupe $ G_k $ op\`ere transitivement dans $
\left(M_{k+1},\omega_{ k+1} \right)\, , $ de
sorte que $ \left(M_k,\omega_ k \right) = \left(M_{k+1},\omega_{ k+1}
\right)\, , $ (voir $\lbrack$NB 1$\rbrack$ Lemme 1.1) On a par
cons\'equent le diagramme commutatif suivant

$$\matrix{\displaystyle H_k & \displaystyle \rightarrow & \displaystyle
H_{k+1} \cr\displaystyle \downarrow\,  & \displaystyle & \displaystyle \downarrow \cr \displaystyle G_k & \displaystyle \rightarrow & \displaystyle
G_{k+1} \cr\displaystyle \downarrow\pi  & \displaystyle  & \displaystyle
\downarrow \cr\displaystyle M_k & \displaystyle \rightarrow & \displaystyle
M_{k+1} \cr}  \leqno (3)$$
Dans le diagramme (3) ci-dessus l'application diff\'erentiable $ \pi $ est la projection canonique de $ G_k $
sur $ G_k/H_k$ et la fl\`eche $ M_k \rightarrow M_{k+1}$ est un diff\'eomorphisme symplectique\, .
Lorsque le sous-groupe $ H_{k+1} $ est inclu dans le sous-groupe $ H_k\, , $
$ H_{k+1} $ est
distingu\'e dans $ H_k $ et sa codimension ibidem vaut 1. La relation entre $
M_k $ et
$ M_{k+1} $ s'exprime en terme de  r\'eduction de phase suivant J. MARSDEN et
WEINSTEIN, $\lbrack$MW$\rbrack$. Plus pr\'ecisement on a le r\'esultat suivant \par
\bigskip
\noindent{\bf Proposition 1.1.1.} {\sl Les notations sont celles du diagramme (2), \/}$ H_{k+1} $
{\sl est inclu dans \/}$ H_k $. {\sl Alors la vari\'et\'e symplectique \/}$
\left(M_k,\omega_ k \right) $ {\sl est une r\'eduite de la
vari\'et\'e symplectique \/}$ \left(M_{k+1},\omega_{ k+1} \right) $ {\sl sous
une action hamiltonienne de \/}$ \bkR $ {\sl dans
\/}$ \left(M_{k+1},\omega_{ k+1} \right) $ \par
\bigskip
\noindent On donnera une d\'emonstration compl\`ete de cette proposition \`a la sous-section . En voici cependant une esquisse. \par
Soit $ \zeta  \in  h_k $ tel que $ H_k $ est produit semi-direct du
sous-groupe \`a un param\`etre $ \{ {\rm expt}  \zeta /t \in  \bkR\} $ par $
H_{k+1}\, . $ Puisque les groupes
en jeu sont tous simplement connexes, l'action de $ \{ {\rm expt}  \zeta /t
\in  \bkR\} $ dans
$ \left(M_{k+1},\omega_{ k+1} \right) $ est hamiltonienne sans point fixe. On
pose
$ \bar  \zeta( x) = {d \over dt}(( {\rm expt}  \zeta)( x))_{\mid t=0}\, . $
\par
\noindent Le sous-groupe $ H_{k+1} $ \'etant connexe, on fixe une primitive $
f_\zeta $ de la 1-forme
ferm\'ee $ i \left(\tilde  \zeta \right)\omega_{ k+1} $ o\`u $ i \left(\tilde
 \zeta \right) $ est le produit int\'erieur par le champ des vecteurs
$ \tilde  \zeta \, . $ La fonction $ f : M_{k+1}\mapsto  \left[h_k/h_{k+1}
\right]^\ast $ d\'efinie par $ (x,\zeta)  \mapsto  f_\zeta( x) $
est une application moment de l'action de $ \{ {\rm expt}  \zeta /t \in
\bkR\} $ dans $ \left(M_{k+1},\omega_{ k+1} \right)\, . $
On verra plus loin que $ \left(M_k,\omega_ k \right) $ est une vari\'et\'e
r\'eduite de cette action
hamiltonienne . \par
\noindent On d\'eduit de la proposition 1.1.1. ci-dessus le diagramme
commutatif
$$ \matrix{\displaystyle H_k & \displaystyle \leftarrow & \displaystyle
H_{k+1} \cr\displaystyle \downarrow & \displaystyle  & \displaystyle
\downarrow \cr\displaystyle G_k & \displaystyle \rightarrow & \displaystyle
G_{k+1} \cr\displaystyle \downarrow & \displaystyle  &
\displaystyle\downarrow \pi \cr\displaystyle{\cal M}_k & \displaystyle
\doubleup{ i \cr \rightarrow \cr} & \displaystyle \left(M_{k+1},\omega_{ k+1}
\right) \cr\displaystyle \pi\downarrow & \displaystyle  & \displaystyle
\cr\displaystyle \left(M_k,\omega_ k \right) & \displaystyle \cdot &
\displaystyle \cr} \leqno (4) $$
dans lequel ${\cal M}_k$ est l'espace homog\`ene $G_k/H_{k+1}$ .\\
Les formes symplectiques $ \omega_ k $ et $ \omega_{ k+1} $ sont li\'ees par
l'\'egalit\'e
$$ i^*w_{k+1} = \pi^*w_k \leqno (5)$$

On voit que dans le cas ci-dessus il n'y a pas de \lq\lq relation
d'incidence\rq\rq\
directe entre $ M_{k+1} $ et $ M_k\, , $ $\lbrack$NB 5$\rbrack$ . On
peut ainsi re\'ecrire (4) sous la
forme
$$ \matrix{\displaystyle H_k & \displaystyle \leftarrow & \displaystyle
H_{k+1} \cr\displaystyle \downarrow & \displaystyle  & \displaystyle
\downarrow \cr\displaystyle G_k & \displaystyle \rightarrow & \displaystyle
G_{k+1} \cr\displaystyle\downarrow \pi & \displaystyle  &
\displaystyle\downarrow \pi \cr\displaystyle M_k & \displaystyle \doubleup{ mw \cr \leftarrow \cr} &
\displaystyle M_{k+1} \cr} \leqno (4^\prime ) $$
o\`u la fl\`eche mw est la r\'eduction symplectique de MARSDEN-WEINSTEIN (voir Proposition 1.1.1).
\par
\medskip
\noindent{\bf D\'efinition 1.1.2.} {\sl Le graphe associ\'e au couple \/}$
(\omega ,F) \in  \Omega^ 2_{\ell}( G) \times  {\cal F}(G) $ {\sl est le
graphe dont l'ensemble \/}$ S $ {\sl des sommets est constitu\'e des \/}$
H_k\, , $ $ G_k\, , $ $ M_k $ {\sl et
dont l'ensemble \/}$ {\cal A} $ {\sl des arr\^etes est constitu\'e des
relations d'inclusion
entre les sommets, }$\lbrack${\sl suivant (1), (2)}$\rbrack${\sl\ et des
projections \/}$ \pi $ $\lbrack${\sl suivant (3) et
(4')}$\rbrack${\sl .\/} \par
\bigskip
Compte tenu des conventions adopt\'ees ci-dessus, la portion du graphe associ\'e \`a $ (\omega ,F)\in  \Omega^
2_{\ell}( G) \times  {\cal F}(G) $ qui
porte les sommets $ H_{k-1}, $ $ H_k $ et $ H_{k+1} $ est une des quatre
configurations
suivantes
$$ (\Irm ) \matrix{\displaystyle H_{k-1} & \displaystyle
\leftarrow & \displaystyle H_k & \displaystyle \leftarrow & \displaystyle H_{k+1}
\cr\displaystyle \downarrow & \displaystyle  & \displaystyle \downarrow &
\displaystyle  & \displaystyle \downarrow \cr\displaystyle G_{k-1} & \displaystyle
\rightarrow & \displaystyle G_k & \displaystyle \rightarrow & \displaystyle G_{k+1}
\cr\displaystyle \downarrow & \displaystyle  & \displaystyle \downarrow &
\displaystyle  & \displaystyle \downarrow \cr \displaystyle M_{k-1}
& \displaystyle \doubleup { mw \cr \leftarrow \cr} & \displaystyle M_k & \displaystyle \doubleup { mw \cr \leftarrow \cr} & \displaystyle M_{k+1}
\cr} \quad ;\quad (\IIrm )\quad \matrix{\displaystyle  H_{k-1} & \displaystyle \rightarrow &
\displaystyle H_k & \displaystyle \rightarrow & \displaystyle  H_{k+1} \cr\displaystyle \downarrow & \displaystyle  & \displaystyle
\downarrow & \displaystyle  & \displaystyle \downarrow \cr\displaystyle G_{k-1} &
\displaystyle \rightarrow & \displaystyle G_k & \displaystyle \rightarrow &
\displaystyle G_{k+1} \cr\displaystyle \downarrow & \displaystyle
& \displaystyle \downarrow  & \displaystyle  & \displaystyle \downarrow
\cr\displaystyle M_{k-1} & \displaystyle \rightarrow & \displaystyle M_k & \displaystyle \rightarrow &
\displaystyle M_{k+1} \cr} \quad ; $$
(6)
$$ (\IIIrm )\quad \matrix{\displaystyle H_{k-1} & \displaystyle \rightarrow & \displaystyle H_k
& \displaystyle \leftarrow & \displaystyle H_{k+1}\cr\displaystyle \downarrow & \displaystyle & \displaystyle \downarrow & \displaystyle  & \displaystyle \downarrow \cr\displaystyle G_{k-1} & \displaystyle \rightarrow & \displaystyle G_k & \displaystyle \rightarrow & \displaystyle G_{k+1} \cr\displaystyle \downarrow & \displaystyle & \displaystyle \downarrow & \displaystyle & \displaystyle \downarrow \cr\displaystyle M_{k+1} & \displaystyle \rightarrow & \displaystyle M_k & \displaystyle \doubleup { mw \cr \leftarrow \cr} & \displaystyle M_k \cr}\quad ;\quad (\Irm \Vrm )\quad \matrix{\displaystyle H_{k-1} & \displaystyle \leftarrow &\displaystyle H_k & \displaystyle \rightarrow & \displaystyle H_{k+1}\cr\displaystyle \downarrow &
\displaystyle  & \displaystyle \downarrow & \displaystyle  &
\displaystyle \downarrow \cr\displaystyle G_{k-1} & \displaystyle \rightarrow & \displaystyle G_k & \displaystyle \rightarrow & \displaystyle G_{k+1} \cr\displaystyle \downarrow &
\displaystyle  & \displaystyle \downarrow & \displaystyle  &
\displaystyle \downarrow \cr\displaystyle M_{k-1} & \displaystyle \doubleup {mw \cr \leftarrow\cr} & \displaystyle M_k & \displaystyle \rightarrow &
\displaystyle M_{k+1}\ . \cr} $$
On va introduire des notions dont la signification g\'eom\'etrique est
donn\'ee dans les sous-sections suivantes. \par
\bigskip
\noindent{\bf D\'efinitions 1.1.3.} (a) {\sl Dans les configurations
(}$\Irm${\sl ) et (}$\IIrm${\sl ) ci-dessus le
sommet \/}$ H_k $ {\sl est dit respectivement \/}r\'egulier non r\'eductible
{\sl et\/} r\'egulier
r\'eductible. {\sl (b) Dans les configurations (}$\IIIrm${\sl ) et (}$\Irm
\Vrm${\sl ) le sommet \/}$ H_k $ {\sl est dit
respectivement\/} singulier attractif {\sl et\/} singulier r\'epulsif
. \par
\bigskip
{\bf Exemple 1.} Soit $\cal G$ l'espace vectoriel de base c,b,a v,u. On munit $\cal G$ du crochet suivant:$ [a,b] = [u,c] = [v,c] = c, [u,a]= a; [v,b]= b;$ les crochets non explicit\'es \'etant nuls. On d\'esigne par ${c^* , b^* , a^* , v^* , u^* } $ la base duale de { c,b,a v,u } et par $F(\cal G) $ le drapeau d\'efini par { c,b,a v,u } et par $ G$ le groupe de Lie connexe et simplement connexe dont l'alg\`ebre de Lie est  $\cal G$.\\
Consid\'erons la 2-forme diff\'erentielle invariante \`a gauche $\omega = a^* \wedge v^* + b^* \wedge u^* + v^* \wedge u^*$; cette forme est ferm\'ee et son noyau est engendr\'e par $c$ . Le sous-espace $vect (a,b,c,v)$ est en fait une sous-alg\`ebre de Lie de  $\cal G$ ; on d\'esigne par $G_4$ le sous-groupe de Lie connexe du groupe $G$ associ\'e \`a $vect (a,b,c,v)$ . La restriction \`a $G_4$ de $\omega$ est la 2-forme $ \omega_4 = a^* \wedge v^* .$ Le noyau de $ \omega_4$ est la sous-alg\`ebre de Lie $vect (c,b)$. On note $G_3$ le sous-groupe de Lie connexe associ\'e \`a la sous-alg\`ebre de Lie $vect (c,b,a)$, la restriction \`a $G_3$ de $\omega$ est nulle. On indexe comme il suit les sous-alg\`ebres de Lie d\'ecrites ci-dessus :\\
$ {\cal H}_5$ = $vect(c)$,\quad ${\cal H}_4$ = $ vect(c, b)$,\quad  ${\cal H}_3$ = $ vect(c,b,a)$ = $ {\cal G}_3 $\\
Le couple $(\omega , F(\cal G))$ donne lieu au diagramme suivant:
$$ \matrix{ \displaystyle H_3 & \displaystyle \leftarrow & \displaystyle H_4 &\displaystyle \leftarrow & \displaystyle H_5\cr\displaystyle \downarrow & \displaystyle & \displaystyle \downarrow & \displaystyle & \displaystyle \downarrow\cr\displaystyle G_3 & \displaystyle \rightarrow & \displaystyle G_4 & \displaystyle \rightarrow & \displaystyle G_5\cr\displaystyle & \displaystyle &  \displaystyle \downarrow & \displaystyle & \displaystyle \downarrow\cr\displaystyle & \displaystyle & \displaystyle M_4 & \displaystyle \doubleup{mw \cr \leftarrow\cr} & \displaystyle M_5 \cr} . $$

{\bf Exemple 2.} On consid\`ere l'espace vectoriel $vect (a,b,c,u)$ muni du crochet suivant: $ [u,c] = -c , [u,b]= b; [u,a]= a;$ les autres crochets sont nuls;${c^* , b^* , a^* , u^* } $ est la base duale de { c,b,a,u }. Soit $\omega = c^* \wedge b^*$ et soit $G$ le groupe de Lie connexe et simplement connexe d'alg\`ebre de Lie ${\cal G} = vect(c,b,a,u); $ $\omega$ d\'efinit dans $G$  un drapeau qui donne lieu au diagramme
$$ H_2 \rightarrow H_3 \rightarrow H_4$$
$$ \matrix {\downarrow & \displaystyle & \downarrow  & \displaystyle & \downarrow\cr} $$
$$ G_2 \rightarrow G_3 \rightarrow G_4$$

{\bf 1.2. DIAGRAMMES POND\'ER\'ES DES 2-FORMES.} \\
Avant de continuer on va fixer quelques conventions. Il est facile
de voir que si deux sommets cons\'ecutifs $ H_k $ et $ H_{k+1} $ sont
r\'eguliers alors
ils sont r\'eguliers de m\^eme type (viz {\sl non r\'eductible\/} ou {\sl
r\'eductible).\/} Cette
remarque conduit \`a contracter les diagrammes en munissant l'arr\^ete joignant deux sommets r\'eguliers cons\'ecutifs d'un indice, e.g. $ \frac{S \qquad S}{(m)}\, ,
$ pour
signifier qu'il y a $ m+1 $ sommets r\'eguliers cons\'ecutifs de type ($S$).
Cela est
commode quand on a l'intention, comme nous nous proposons de faire,
d'associer \`a chaque graphe orient\'e l'analogue de diagramme de Dynkin pour
les syst\`emes des racines, $\lbrack$B$\rbrack$. \par
Soit $ (\omega ,F)\in  \Omega^ 2_{\ell}( G) \times  {\cal F}(G)\,$ . Notons $ {\rm gr}(\omega ,F) $ le
graphe orient\'e associ\'e au couple $ (\omega ,F)\, . $ On consid\`ere
maintenant le diagramme $ d(w,F)$ dont les sommets $ S_k $ ont pour coordonn\'ees les couples $
\left(G_k,H_k \right)\, . $ On convient de relier les deux sommets cons\'ecutifs $ S_k $ et $ S_{k+1} $ par deux arr\^etes orient\'ees exprimant les relations d'inclusion ad hoc; on notera que dans les diagrammes commutatifs (5) et (6) la fl\`eche joignant $H_k$ \`a $H_{k-1}$ (respectivement \`a $H_{k+1}$) et celle joignant $M_k$ \`a $M_{k-1}$ (respectivement $M_{k+1}$) ont la m\^eme orientation. Ainsi compte t\'enu de notre convention de contraction des diagrammes les configurations ($\Irm$), ($\IIrm$), ($\IIIrm$) et ($\Irm \Vrm$) donnent lieu aux quatre diagrammes suivants
$$ \matrix{\displaystyle (\Irm ) & \displaystyle \doubleup{ S_{k-1} \cr {\rm
O} \cr}^{ \leftarrow}_{ \rightarrow} \doubleup{ S_k \cr {\rm O} \cr}^{
\leftarrow}_{ \rightarrow} \doubleup{ S_{k+1} \cr {\rm O} \cr} \hfill&
\displaystyle ;\quad (\IIrm ) \hfill& \displaystyle {\rm O}^{\rightarrow}_{
\rightarrow} \doubleup{ S_k \cr {\rm O} \cr}^{ \rightarrow}_{ \rightarrow}
{\rm O} & \displaystyle ; \cr\displaystyle  & \displaystyle  & \displaystyle
& \displaystyle  & \displaystyle \cr\displaystyle (\IIIrm ) & \displaystyle
{\rm O}^{\rightarrow}_{ \rightarrow} \doubleup{ S_k \cr {\rm O} \cr}^{
\leftarrow}_{ \rightarrow} {\rm O}  \hfill& \displaystyle ;\quad (\Irm \Vrm )
\hfill& \displaystyle {\rm O}^{\leftarrow}_{ \rightarrow} \doubleup{ S_k \cr
{\rm O} \cr}^{ \rightarrow}_{ \rightarrow} {\rm O} & \displaystyle . \cr} $$
On peut convenir de joindre $ S_{k-1} $ \`a $ S_k $ par une ou deux arr\^etes
suivant que
l'on a $ H_{k-1}\subset  H_k $ ou $ H_k\subset  H_{k-1}\, , $ de sorte que les
diagrammes ci-dessus
deviennent les suivantes
$$ \matrix{\displaystyle (\Irm ) & \displaystyle {\rm O}^{\leftarrow}_{
\rightarrow} \doubleup{ S_k \cr {\rm O} \cr}^{ \leftarrow}_{ \rightarrow}
{\rm O}   & \displaystyle ;\quad (\IIrm ) \hfill& \displaystyle {\rm O}
\rightarrow  \doubleup{ S_k \cr {\rm O} \cr} \rightarrow {\rm O}   &
\displaystyle ; \cr\displaystyle  & \displaystyle  & \displaystyle  &
\displaystyle  & \displaystyle \cr\displaystyle (\IIIrm ) & \displaystyle
{\rm O} \rightarrow \doubleup{ S_k \cr {\rm O} \cr}^{ \leftarrow}_{
\rightarrow} {\rm O}   & \displaystyle ;\quad (\Irm \Vrm ) \hfill&
\displaystyle {\rm O}^{\leftarrow}_{ \rightarrow} \doubleup{ S_k \cr {\rm O}
\cr} \rightarrow {\rm O}   & \displaystyle . \cr} \leqno (7) $$
On voit que les sommets $ S_k $ des configurations ($\Irm$) et ($\IIrm$) sont
r\'eguliers
dans le sens \'evident qu'il y a autant d'arr\^etes d'origine $ S_k $ que
celles
d'extr\'emit\'es $ S_k\, . $ En contraste il y a plus d'arr\^etes d'origine $
S_k $ (resp.
d'extr\'emit\'e $ S_k) $ que celles d'extr\'emit\'e $ S_k $ (resp. d'origine $
S_k) $ dans le cas
($\IIIrm$) (resp. ($\Irm \Vrm$)). \par
\bigskip
\noindent{\bf Poids des sommets.} On va attacher au sommet $ S_k=
\left(H_k,G_k \right) $ le poids
$ p_k= p \left(S_k \right) = {\dim H_k \over {\rm codim} \ H_k+1} = {\dim H_k
\over \dim M_k+1}\, . $ La notion de poids est li\'ee \`a des
propri\'et\'es g\'eom\'etriques de $ (\omega ,F)\, ; $ on y reviendra dans les
num\'eros
suivants. \par
\bigskip
{\bf 1.3. Connexit\'e et simplicit\'e des diagrammes pond\'er\'es.} \\
Soit $ \omega  \in  \Omega^ 2_{\ell}( G)\, . $ Soit $ g $
l'alg\`ebre de Lie du groupe de Lie $ G $ et soit $ h $ la sous-alg\`ebre de
Lie
de $ g $ qui engendre le noyau de la 2-forme diff\'erentielle $ \omega \, . $
Si le
groupe de Lie $ G $ est nilpotent, alors il existe dans $ G $ un sous-groupe
de
Lie $ {\cal G} $ de dimension $ = \dim H + {1 \over 2}\ {\rm codim} \ H $
jouissant des deux propri\'et\'es
suivantes
$$ ( {\rm a})\quad H \subset  {\cal G}\, ; ( {\rm b})\quad i^\ast \omega  = 0
$$
o\`u $ i $ est l'homomorphisme inclusion de $ {\cal G} $ dans $ G\, ; $ le th\'eor\`me 1.4.1 donnera entre autres un analogue r\'esoluble de cette propri\'et\'e.
Naturellement par passage
\`a la vari\'et\'e quotient $ M = G/H $ on obtient la vari\'et\'e symplectique
not\'ee $ (M,\omega) $
et la sous-vari\'et\'e $ L = {\cal G}/H $ co\"\i ncide avec la feuille passant
par $ x_o= \pi( H) $
d'un feuilletage lagrangien $ {\cal F} $ invariant par l'action de $ G $ dans
$ G/H\, . $ Il
existe une suite de composition $ F \in  {\cal F}(G) $ contenant $ H $ et $
{\cal G}\, , $ de sorte que
le diagramme pond\'er\'e $ d(\omega ,F) $ associ\'e \`a $ gr(\omega ,F) $ a
une des formes suivantes
$$ \matrix { \displaystyle & \displaystyle & \displaystyle S & \displaystyle &\displaystyle S^{\prime} & \displaystyle & \displaystyle \cr \displaystyle O & \displaystyle \rightarrow & \displaystyle O & \displaystyle  ^{\leftarrow}_{\rightarrow} & \displaystyle O & \displaystyle \rightarrow & \displaystyle O\cr\displaystyle & \displaystyle (m) & \displaystyle & \displaystyle (n) & \displaystyle & \displaystyle & \displaystyle & \cr}\leqno (i)  $$
$$ \matrix { \displaystyle & \displaystyle & \displaystyle S & \displaystyle & \displaystyle \cr\displaystyle O & \displaystyle \rightarrow &  \displaystyle O & \displaystyle ^{\leftarrow}_{\rightarrow} & \displaystyle O\cr\displaystyle & \displaystyle (m) & \displaystyle & \displaystyle (n) & \displaystyle \cr } \leqno (ii) $$

Les pr\'esentations (i) et (ii) ci-dessus ayant les signification convenues
au num\'ero 1.2 en d'autres termes l'indice (m) (respectivement (n)) signifie qu'on a m+1 sommets cons\'ecutifs li\'es comme indiqu\'e ci-dessus. On peut r\'esumer (i) et (ii) sous la forme de \par
\bigskip
\noindent{\bf Proposition 1.3.1.} {\sl Soit \/}$ G $ {\sl un groupe de Lie
nilpotent connexe, pour toute
2-forme ferm\'ee \/}$ \omega  \in  \Omega^ 2_{\ell}( G)\, , $ {\sl il existe
\/}$ F \in  {\cal F}(G) $ {\sl tel que \/}$ d(\omega ,F) $ {\sl poss\`ede au
plus deux sommets singuliers. Lorsque \/}$ d(\omega ,F) $ {\sl poss\`ede un
seul sommet singulier, celui-ci est attractif }. \par
\bigskip
\noindent Cette proposition 1.3.1 permet une bonne compr\'ehension des graphes
des
2-formes ferm\'ees $ \omega  \in  \omega^ 2_{\ell}( G) $ quand $ G $ est
compl\`etement r\'esoluble. \par
\noindent Il faut pour cela reformuler le lemme 1.1.1 de $\lbrack$NB
1$\rbrack$ de la fa\c con qui
suit. \par
\smallskip
\noindent{\bf Proposition 1.3.2.} {\sl Soit \/}$ G $ {\sl un groupe de Lie
compl\`etement r\'esoluble non
nilpotent. Soit \/}$ \omega  \in  \Omega^ 2_{\ell}( G) $ {\sl une 2-forme
 non nulle. On note \/}$ H $ {\sl le
sous-groupe de Lie connexe dont l'alg\`ebre de Lie \/}$ h $ {\sl engendre le
noyau
de \/}$ \omega \, . $ {\sl Si \/}$ h $ {\sl est non nulle, alors les
affirmations suivantes sont
\'equivalentes.\/} \par
\noindent$ \left(a_1 \right) $ {\sl Le seul sous-groupe distingu\'e de \/}$ G
$ {\sl qui contient \/}$ H $ {\sl est \/}$ G $ {\sl lui-m\^eme.\/} \par
\noindent$ \left(a_2 \right) $ {\sl Le sous-groupe des commutations \/}$
{\cal D}G $ {\sl op\`ere transitivement dans \/}$ G/H  $ \par
\bigskip
\noindent{\bf Esquisse de d\'emonstration.} L'\'equivalence de $ \left(a_1
\right) $ et $ \left(a_2 \right) $ se d\'eduit
ais\'ement de la commutativit\'e du diagramme suivant
$$ \matrix{\displaystyle H \cap  {\cal D}G & \displaystyle \hookrightarrow &
\displaystyle{\cal D}F & \displaystyle \rightarrow & \displaystyle{\cal D}G/H
\cap  {\cal D}G \cr\displaystyle \downarrow & \displaystyle  & \displaystyle
\downarrow & \displaystyle  & \displaystyle \downarrow \cr\displaystyle H &
\displaystyle \hookrightarrow & \displaystyle G & \displaystyle \rightarrow &
\displaystyle G/H \cr\displaystyle \downarrow & \displaystyle  &
\displaystyle \downarrow & \displaystyle  & \displaystyle \downarrow
\cr\displaystyle H/H\cap{\cal D}G & \displaystyle \hookrightarrow &
\displaystyle G/{\cal D}G & \displaystyle \rightarrow & \displaystyle{\cal A}
\cr} $$
dont les trois colonnes ainsi que les trois lignes sont les fibrations
diff\'erentiables . \par
\smallskip
\noindent On se propose avant d'avancer de tirer des propositions 1.3.1 et
1.3.2 des
informations concernant le nombre et les natures des sommets singuliers
des diagrammes pond\'er\'es $ d(\omega ,F) $ quand $ G $ est compl\`etement
r\'esoluble non
nilpotent. En effet, il r\'esulte de la proposition 1.3.2. que pour toute
2-forme  $ \omega  \in  \Omega^ 2_{\ell}( G) $ il existe une
filtration $ F $ de $ G $ par des
sous-groupes $ G_k $ telle que le diagramme pond\'er\'e $ d(\omega ,F) $ a une
des deux
formes qui suivent
$$ \matrix{\displaystyle {\rm O} \rightarrow {\rm O}^{\leftarrow}_{
\rightarrow} {\rm O} \rightarrow {\rm O}^{\leftarrow}_{ \rightarrow} {\rm O}
\rightarrow {\rm O} \ , \cr\displaystyle \cr\displaystyle {\rm O} \rightarrow
{\rm O}^{\leftarrow}_{ \rightarrow} {\rm O} \rightarrow {\rm
O}^{\leftarrow}_{ \rightarrow} {\rm O\ .} \cr} $$
On verra que des difficult\'es de compr\'ehension des dynamiques symplectiques
de la vari\'et\'e symplectique $ (M,\omega) $ d\'eduite de $ \omega  \in
\Omega^ 2_{\ell}( G)\, , $ certaines sont
li\'ees \`a la nature et au nombre des sommets singuliers des diagrammes
pond\'er\'es  $ d(\omega ,F)\, . $ Il en est ainsi par exemple quand on est
int\'eress\'e par
l'existence des feuilletages lagrangiens invariants par un groupe des
diff\'eomorphismes symplectiques de $ (M,\omega) \, . $ \par
\noindent On observe enfin qu'un sommet singulier de poids nul est
n\'ecessairement r\'epulsif . Il va de soit qu'un sommet de poids nul est singulier r\'epulsif . \par
\smallskip
Les consid\'erations et observations ci-dessus motivent les
d\'efinitions qui suivent. \par
\bigskip
\noindent{\bf D\'efinitions 1.3.3.} {\sl (a) Un diagramme pond\'er\'e \/}$
d(\omega ,F) $ {\sl est \/}dit connexe{\sl\  s'il ne poss\`ede pas de sommet singulier de poids nul }.
\par
\noindent{\sl (b) Un diagramme pond\'er\'e \/}$ d(\omega ,F) $ est
semi-nilpotent{\sl\  si ses sommets singuliers de poids nul \/}sont nilpotents . \par
\bigskip
\noindent De la d\'efinition 1.3.3. r\'esulte que tout diagramme pond\'er\'e $
d(\omega ,F) $ est
r\'eunion de ses composantes connexes. Naturellement deux composantes connexes distinctes ont au plus un sommet commun qui est de poid  nul . Par analogie
\`a la  th\'eorie des  racines des groupes de Lie semi-simples complexes on va poser. \par
\bigskip
\noindent{\bf D\'efinition 1.3.4.} {\sl Un diagramme pond\'er\'e \/}$
d(\omega ,F) $ {\sl est\/} dit simple {\sl s'il est connexe et ne poss\`ede pas de sommet r\'epulsif }. \par
\bigskip
\noindent Il revient au m\^eme de dire que $ d(\omega ,F) $ est simple  s'il
poss\`ede un seul sommet singulier. Ce dernier est alors n\'ecessairement
attractif. Un tel diagramme simple a donc la forme suivante
.
$$ {\rm O} \rightarrow {\rm O}^{\leftarrow}_{ \rightarrow} {\rm O} \quad .  $$
Pour illustrer les notions qui viennent d'\^etre introduites \`a l'aide de  quelques exemples, la d\'efinition qui suit sera utile \par
\bigskip
\noindent{\bf D\'efinitions 1.3.5.} {\sl (a) Un diagramme pond\'er\'e \/}$
d(\omega ,F) $ {\sl dont les sommets
singuliers de poids nul sont distingu\'es dans \/}$ G $ {\sl est dit
\/}semi-normal. \par
\noindent{\sl (b) Un diagramme pond\'er\'e semi-normal dont les compossantes
connexes sont simples est dit \/}semi-simple . \par
\bigskip
\noindent Un diagramme pond\'er\'e semi-simple a la forme suivante

$$ \matrix { \displaystyle O & \displaystyle \rightarrow &  \displaystyle O & \displaystyle ^{\leftarrow}_{\rightarrow} & \displaystyle O & \displaystyle \cdots & \displaystyle O & \displaystyle \rightarrow &  \displaystyle O & \displaystyle ^{\leftarrow}_{\rightarrow} & \displaystyle O\cr\displaystyle & \displaystyle (m_1) & \displaystyle & \displaystyle (n_1) & \displaystyle & \displaystyle & \displaystyle & \displaystyle (m_{\ell}) & \displaystyle & \displaystyle (n_{\ell}) & \displaystyle\cr }  $$

o\`u les pointill\'es repr\'esentemt les diagrammes simples pond\'er\'es $$ \matrix { \displaystyle O &\displaystyle \rightarrow & \displaystyle O & \displaystyle ^{\leftarrow}_{\rightarrow} & \displaystyle O \cr } .$$

Bien entendu les sommets singuliers r\'epulsifs sont de poids nul . \par
\bigskip
{1.4. {\bf Des exemples.} \\
\smallskip
\noindent{\bf Exemple 1.} On note $ G = \left\{ A = \left[
\matrix{\displaystyle \alpha & \displaystyle u & \displaystyle v
\cr\displaystyle 0 & \displaystyle \alpha^ 2 & \displaystyle w
\cr\displaystyle 0 & \displaystyle 0 & \displaystyle \alpha^ 3 \cr} \right],
(\alpha ,u,v,w)\in  \bkR^{ >0}\times  \bkR^ 3 \right\} \, . $ \par
\noindent On note $ A(\alpha ,u,v,w) $ les \'el\'ements de $ G\, . $ \par
\noindent Soit  $\omega = 2d\alpha \wedge dv - du \wedge dw $.
 On va filtrer $ G $ par les sous-groupes
$$ G_1= \{ A(\alpha ,u,v,w) / \alpha = 1 , u = w = 0\}  ; G_2= \{ A(\alpha ,u,v,w) / \alpha = 1 , u =0  \}  ; G_3= \{ A(\alpha ,u,v,w) / \alpha = 1 \} \, .$$
On voit ais\'ement que
$$ H_1= G_1, H_2= G_2\ {\rm et} \ H_3= \{ A(1,u,0,0)\} \ . $$
Le diagramme associ\'e \`a la filtration ci-dessus est
$$ {\rm O} \doublelow{ \rightarrow \cr (2) \cr} {\rm O} \doublelow{^{
\leftarrow}_{ \rightarrow} \cr (1) \cr} {\rm O\ .} $$
C'est un diagramme simple. \par
\bigskip
\noindent{\bf Exemple 2.} On conserve $ G $ et $ \omega $ de l'exemple 1. On
filtre $ G $ par les
sous-groupes suivants
$$ G_1= \{ A(\alpha ,u,v,w) / u = v = w = 0\}  ; G_2= \{ A(\alpha ,u,v,w) / u = v = 0 \}  ; G_3= \{ A(\alpha ,u,v,w) / u = 0 \} {\rm \ .} $$
On a alors $ H_1= G_1; H_2= \{ A(1,0,0,0)\}  ; H_3= \{ A(1,0,v,0)\} \, . $
\par
\noindent Le diagramme correspondant \`a cette filtration est
$$ {\rm O} \rightarrow {\rm O}^{\leftarrow}_{ \rightarrow} \doubleup{ S_2 \cr
{\rm O} \cr} \rightarrow {\rm O}^{\leftarrow}_{ \rightarrow} {\rm O} \ .  $$
Le sommet $ S_2= \left(H_2,G_2 \right) $ est de poids nul ; $ d(\omega ,F) $
n'est donc pas connexe. Il
a deux composantes simples et le sommet singulier $ S_2 $ est distingu\'e ; $
d(\omega ,F) $
est semi-simple. \par
\bigskip
\noindent{\bf Exemple 3.} Soit $ G $ le groupe de Lie connexe des matrices r\'eelles $ 5
\times  5 $ dont l'alg\`ebre de
Lie $ g $ est l'espace des matrices
$$ X(z,y,x,u,v) = \left[ \matrix{\displaystyle 2u & \displaystyle x+v &
\displaystyle -y & \displaystyle -2z & \displaystyle -y \cr\displaystyle 0 &
\displaystyle u & \displaystyle v & \displaystyle -y & \displaystyle -x
\cr\displaystyle 0 & \displaystyle 0 & \displaystyle u & \displaystyle -x-v &
\displaystyle u \cr\displaystyle 0 & \displaystyle 0 & \displaystyle 0 &
\displaystyle 0 & \displaystyle 0 \cr\displaystyle 0 & \displaystyle 0 &
\displaystyle 0 & \displaystyle 0 & \displaystyle 0 \cr} \right]\ . $$
Notons $ \left(e_5,e_4,e_3,e_2,e_1 \right) $ la base canoniquement associ\'ee
aux coordonn\'ees  $ (z,y,x,u,v) $ et $ \left(\theta_ 5,\theta_ 4,\theta_ 3,\theta_ 2,\theta_ 1 \right) $ la base duale. On consid\`ere la 2-forme
$ \omega  = d\theta_ 5\, . $ On pose $ X(z,y,x,u,v) = ze_5 + ye_4 + xe_3 + ue_2 + ve_1 $ . On a alors
$$ \omega( X(z,y,x,u,v),X(z^\prime ,y^\prime ,x^\prime ,u^\prime ,v^\prime))
= xy^\prime  - x^\prime y + 2(uz^\prime - u^\prime z) + vy^\prime  - v^\prime
y\, . $$
Le noyau de $ \omega $ est engendr\'e par la sous-alg\`ebre de $ g $ d\'efinie
par
$$ z = y = u = 0\, , x + v = 0\, . $$
On va filtrer l'alg\`ebre de Lie $ g $ du groupe $ G $ par les sous-alg\`ebres suivantes
\medskip
$$ F :\ =  g(0,0,x,0,-x) \subset  g(z,0,x,0,-x) \subset  g(z,0,x,0,v) \subset
g(z,y,x,0,v) \subset  g\, . $$
Les notations ci-dessus ne pr\'etant pas \`a ambiguit\'e, puisque
$$ g(0,0,x,0,-x) = \{ X(0,0,x,0,-x)/ x \in \bkR \} \, . $$
Le noyau de $ \omega $ est engendr\'e par $ g(0,0,x,0,-x)\, . $ On voit que le
noyau de la restriction $ \omega( z,y,x,0,v) = \omega_{ \mid g(z,y,x,x,0,v)} $ est $ g(z,0,x,0,-x)\, . $ \par
\noindent Le diagramme $ d(\omega ,F) $ est simple et a la forme suivante
$$ {\rm O}  \doublelow{ \rightarrow \cr (2) \cr} \doubleup{ S_5 \cr {\rm O}
\cr}^{ \leftarrow}_{ \rightarrow} {\rm O} \, . $$
On conserve $ G $ et $ \omega \,  $ comme ci-dessus, puis on filtre $ g $ par les sous-alg\`ebres comme ci-dessous
\medskip
$$ F :\ = g(z,0,0,0,0) \subset  g(z,y,0,0,0) \subset  g(z,y,x,0,0) \subset
g(z,y,x,u,0) \subset  g\, . $$
La restriction \`a $ g(z,y,x,u,0) $ de $ \omega $ est alors r\'eguli\`ere et
le noyau de
$ \omega( z,y,x,0,0) $ est $ g(z,0,0,0,0)\, . $ Le diagramme $ d(\omega ,F) $
est
$$ {\rm O} \rightarrow {\rm O}^{\leftarrow}_{ \rightarrow} \doubleup{ S_4 \cr
{\rm O} \cr} \rightarrow {\rm O} \ . $$
On peut avoir un diagramme connexe non semi-simple en filtrant $ g $
comme ci-dessous
\medskip
$$ F :\ = g(0,y,0,u,0) \subset  g(z,y,0,0,0) \subset  g(z,y,0,u,0) \subset
g(z,y,x,u,-x) \subset  g\, . $$
Le noyau de $ \omega( z,u,x,u,-x) $ est $ g(0,y,x,0,-x)\, ; $ le noyau de $
\omega( z,y,0,u,0) $
est $ g(0,y,0,0,0)\, . $ On a visiblement $ \omega( z,y,0,0,0) = 0\, . $ Le
diagramme $ d(\omega ,F)$ est le suivant
$$ {\rm O} \rightarrow {\rm O}^{\leftarrow}_{ \rightarrow} {\rm O}
\rightarrow {\rm O}^{\leftarrow}_{ \rightarrow} {\rm O\, ,} $$
Ce diagramme est connexe mais non semi-simple . \par
\smallskip
\noindent On voit qu'\`a divers choix de $ F \in  {\cal F}(G) $ correspondent
des diagrammes de
types diff\'erents. On donnera des significations g\'eom\'etriques aux
configurations int\'eressantes. \par
\noindent Nous allons \'etablir un r\'esultat qui joue un r\^ole utile dans la
dynamique
des vari\'ets K\"ahl\'eriennes compactes. \par
\bigskip
\noindent \overfullrule=0mm \par
\noindent{\bf 1.4. R\'eduction des diagrammes semi-nilpotents.} \par
Soit $ G $ un groupe de Lie compl\`etement r\'esoluble. Notre objectif est
la recherche des possibilit\'es d'avoir pour une forme ferm\'ee $ \omega \in
\Omega^ 2_{\ell}( G) $ un diagramme simple. \par
\bigskip
\noindent{\bf Th\'eor\`eme 1.4.1.} {\sl Soit \/}$ G $ {\sl un groupe de Lie
compl\`etement r\'esoluble. Soit
\/}$ \omega  \in  \Omega^ 2_{\ell}( G) $ . {\sl Si
\/}$ \omega $ {\sl poss\`ede un diagramme pond\'er\'e semi-simple
et semi-nilpotent \/}$ d(\omega ,F)\, , $ {\sl alors on peut \lq\lq
d\'eformer\rq\rq\ \/}$ F $ {\sl en une filtration \/}$ F_0 $
{\sl tel que \/}$ d \left(\omega ,F_0 \right) $ {\sl est un diagramme
simple\/}
$$ {\rm O} \rightarrow {\rm O}^{\leftarrow}_{ \rightarrow} {\rm O} \quad. $$
{\sl D\'emonstration.\/} D'apr\`es l'hypoth\`ese $ \omega $ poss\`ede un
diagramme semi-simple $ d(\omega ,F) $ de la forme

$$ \matrix { \displaystyle O &\displaystyle\rightarrow &\displaystyle O & \displaystyle^{\leftarrow}_{\rightarrow} & \displaystyle O &\displaystyle \rightarrow &  \displaystyle O & \displaystyle ^{\leftarrow}_{\rightarrow} & \displaystyle O & \displaystyle \cdots & \displaystyle O & \displaystyle \rightarrow &  \displaystyle O & \displaystyle ^{\leftarrow}_{\rightarrow} & \displaystyle O\cr\displaystyle & \displaystyle (m_1) & \displaystyle & \displaystyle (n_1) & \displaystyle & \displaystyle (m_2) & \displaystyle & \displaystyle (n_2) & \displaystyle & \displaystyle & \displaystyle & \displaystyle (m_{\ell}) & \displaystyle & \displaystyle (n_{\ell}) & \displaystyle\cr }  $$

Les sommets singuliers r\'epulsifs \'etant distingu\'es et nilpotents. Soit $
2n $
le rang de $ \omega $ et soit $ 2n+\ell $ la dimension de $ G\, . $ Soit $ g $
l'alg\`ebre de
Lie du groupe de Lie $ G\, . $ Notons comme auparavant $ h $ la sous-alg\`ebre
de Lie de $ g $ qui engendre le noyau de $ \omega \, , $ la dimension de $ h $
est
\'egale \`a $ \ell \, . $ Soit $ H $ le sous-groupe de Lie connexe associ\'e
\`a $ h\, . $ \par
\noindent Il s'agit de montrer qu'il existe une filtration de $ G $ par les
sous-groupes $ G_k $ dont celui de dimension $ n + \ell \, , $ viz $ G_{n+\ell} \, , $
jouit des deux propri\'et\'es suivantes
$$ H \subset G_{n+\ell} ; \leqno (p1)$$
$$ i^*_{n+\ell}\omega = 0 . \leqno (p2)$$
On part du diagramme semi-simple semi-nilpotent ci-dessus, soit

$$ \matrix { \displaystyle O & \displaystyle \rightarrow &  \displaystyle O & \displaystyle ^{\leftarrow}_{\rightarrow} & \displaystyle O & \displaystyle \cdots & \displaystyle O & \displaystyle \rightarrow &  \displaystyle O & \displaystyle ^{\leftarrow}_{\rightarrow} & \displaystyle O\cr\displaystyle & \displaystyle (m_1) & \displaystyle & \displaystyle (n_1) & \displaystyle & \displaystyle & \displaystyle & \displaystyle (m_{\ell}) & \displaystyle & \displaystyle (n_{\ell}) & \displaystyle\cr }  $$

On peut supposer sans perte de g\'en\'eralit\'e que $ d(\omega ,F) $ a
seulement deux composantes connexes simples. Cette restriction est l\'egitime en vertu des propositions 1.3.1 et 1.3.2. Le diagramme pond\'er\'e $ d(\omega ,F) $ est d\'esormais suppos\'e \^etre le suivant

$$ \matrix { \displaystyle & \displaystyle & \displaystyle S^\prime & \displaystyle & \displaystyle S & \displaystyle & \displaystyle S^{''} & \displaystyle & \displaystyle\cr\displaystyle O & \displaystyle \rightarrow &  \displaystyle O & \displaystyle ^{\leftarrow}_{\rightarrow} & \displaystyle O & \displaystyle \rightarrow & \displaystyle O & \displaystyle^{\leftarrow}_{\rightarrow} & \displaystyle O\cr\displaystyle & \displaystyle (m_1) & \displaystyle & \displaystyle (n_1) & \displaystyle & \displaystyle (m_2) & \displaystyle & \displaystyle (n_2) \cr } $$
dans lequel l'unique sommet singulier r\'epulsif $ S $ est de poids nul ; il a pour coordonn\'ees le couple $ \left(G_{2m}, \{ e\} \right) $ o\`u $ G_{2m} $ est un sous-groupe distingu\'e nilpotent du groupe $ G $, $ e $ \'etant l'\'el\'ement neutre de $ G\, . $ \par
\noindent Naturellement la 2-forme $ \omega_{ 2m}= i^\ast_{ 2m}\omega $ est
symplectique dans $ G_{2m}\, . $ Il en
r\'esulte que $ G $ est produit semi-direct $ G = G_{2m}$ $\times\hskip-0,08cm\vbox{\hbox{$^{_|}$}\vskip-0,045cm}$ $ G^o $ o\`u $
G^o $ est un sous-groupe
de $ G $ de dimension $ 2(n-m)+\ell \, ; $ $ G^o $ est l'\lq\lq
orthogonal\rq\rq\ de $ G_{2m} $ relativement \`a $ \omega \, . $
Il resulte de la d\'ecomposition en produit semi-direct
$$ G = G_{2m}\hbox{$\times\hskip-0,08cm\vbox{\hbox{$^{_|}$}\vskip-0,045cm}$} G^o $$
que la forme $ \omega $ se d\'ecompose en deux termes
$$ \omega  = \omega_{ 2m}+ \omega^ o $$
o\`u $ \omega^ o $ est la restriciton \`a $ G^o $ de $ \omega \, $ et $ \omega_{2m}$ est la restriction \`a $ G_{2m}$ de $ \omega $\par
\noindent Le sommet singulier $ S^{\prime\prime} $ correspond (pour la
filtration $ F $ en jeu) \`a un
sous-groupe $ G_{n+2m+\ell} $ dans $ G\, . $ Du point de vue du produit
semi-direct
$ G = G_{2m}$ $\times\hskip-0,08cm\vbox{\hbox{$^{_|}$}\vskip-0,045cm}$ $G^o\, , $ le sous-groupe $ G_{n+2m+\ell} $ a la forme
$$ G_{n+2m+\ell} = G_{2m}\hbox{$\times\hskip-0,08cm\vbox{\hbox{$^{_|}$}\vskip-0,045cm}$}  G^o_{n-m+\ell} $$
o\`u le sous-groupe $ G^o_{n-m+\ell} $ de $ G^o $ satisfait les deux
conditions suivantes
$$ H \subset G^{\circ}_{n-m+{\ell}} ; \leqno (c1)$$
$$ i^*_{n-m+{\ell}}\omega^{\circ} = 0 . \leqno (c2)$$

Nous sommes maintenant en mesure de construire une filtration $ F_o $ telle
que
$ d \left(\omega ,F_o \right) $ soit simple. Il suffit pour cela de
d\'emontrer le Lemme qui suit. \par
\bigskip
\noindent{\bf Lemme 1.4.2.} {\sl Le sous-groupe \/}$ G_{2m} $ {\sl poss\`ede
une filtration \/}$ \tilde  F \left(G_{2m} \right) = G_1\subset  G_k\cdot
\cdot \cdot $
{\sl jouissant des propri\'et\'es suivantes\/} \par
\noindent{\sl (i) chaque \/}$ G_k $ {\sl est stable par les automorphismes
int\'erieurs de \/}$ G $ {\sl d\'etermin\'es
par les \'el\'ements de \/}$ G^o\, . $ \par
\noindent{\sl (ii) \/}$ i^\ast_ m\omega_{ 2m}= 0\, , $ {\sl o\`u \/}$ i_m $
{\sl est l'homomorphisme inclusion de \/}$ G_m $ {\sl dans \/}$ G_{2m}\, . $
\par
\bigskip
\noindent{\sl D\'emonstration du Lemme 1.4.2.\/} Soit $ g_m $ l'alg\`ebre de
Lie du sous-groupe $ G_{2m} $
et soit $ g^o $ celle du sous-groupe de Lie $ G^o\, . $ L'alg\`ebre de Lie $ g
$ de $ G $ est
produit semi-direct de $ g^o $ par $ g_{2m}\, . $ Notons $ a $ la
sous-alg\`ebre de Lie de $ g^o $
qui correspond au sous-groupe de Lie $ G^o_{n-m+\ell} \, . $ \par
\bigskip
\noindent{\bf 1\`ere Etape. }On observe que le noyau de la restriction de $
\omega $ \`a $ g_{2m}\times  a $ est $ a\, . $ \par
\noindent On fixe un id\'eal de codimension 1, $ g_{2m-1}\, , $ dans $ g_{2m}
$ astreint \`a v\'erifier la
condition $ \left[g^o, g_{2m-1} \right] \subset  g_{2m-1}\, . $ \par
\noindent Un tel id\'eal existe parceque $ G $ est compl\`etement r\'esoluble
d'une part,
et d'autre part toute sous-alg\`ebre de codimension 1 dans $ g_{2m} $ est un
id\'eal de $ g_{2m}\, . $ \par
\noindent Notons $ \omega_{2m-1} $ la restriction \`a $ g_{2m-1} $ de $
\omega_{ 2m}\, , $ le noyau de $ \omega_{ 2m-1} $ est de
dimension 1, notons $ h_1 $ ce noyau. Puisque $ g^o $ agit dans $
\left(g_{2m},\omega_{ 2m} \right) $ comme
sous-alg\`ebre de $ sp \left(\omega_{ 2m} \right) $ on a n\'ecessairement
$$ \left[g^o,h_1 \right] \subset  h_1\, . $$
Naturellement le noyau de la restriction \`a $ g_{2m-1}\times  a $ de $
\omega $ est la
sous-alg\`ebre $ h_1\times  a\, . $ Il d\'ecoule de l'inclusion
$$ \left[g^o,h_1 \right] \subset  h_1 $$
que l'on a aussi
$$ \left[g^o,nor(h_1) \right] \subset  nor(h_1) \leqno (8) $$
o\`u $ {\rm nor} \left(h_1 \right) $ est le normalisateur de $ h_1 $ dans $
g_{2m-1}\, . $ \par
\noindent On d\'eduit de (8) l'existence d'une sous-alg\`ebre de codimension 1
dans $ g_{2m-1}\, , $
soit $ g_{2m-2}\, , $ qui v\'erifie les conditions
$$ h_1\subset  g_{2m-2}\ {\rm et} \ \left[g^o,g_{2m-2} \right] \subset
g_{2m-2}\, . $$
L'existence de $ g_{2m-2} $ r\'esulte entre autres de la nilpotence de $
g_{2m-1}\, . $ \par
\noindent Le noyau de la restriction \`a $ g_{2m-2} $ de $ \omega_{ 2m-1} $
est de dimension 2 ; notons ce
noyau $ h_2\, . $ On a comme ci-dessus les inclusions
$$ \matrix{\displaystyle \left[g^o,h_2 \right] \subset  h_2, \cr\displaystyle
\left[g^o, {\rm nor} \left(h_2 \right) \right] \subset  {\rm nor} \left(h_2
\right) \cr} $$
o\`u $ {\rm nor} \left(h_2 \right) $ est le normalisateur de $ h_2 $ dans $
g_{2m-2}\, . $ \par
\noindent Le noyau de la restriction \`a $ g_{2m-2}\times  a $ de $ \omega $
est la sous-alg\`ebre $ h_2\times  a\, . $
Bien \'evidemment on a de nouveau l'existence d'une sous-alg\`ebre de Lie de
codimension 1 dans $ g_{2m-2}\, , $ soit $ g_{2m-3} $ qui contient $ h_2 $ et
v\'erifie
$$ \left[g^o,g_{2m-3} \right] \subset  g_{2m-3} $$
On it\`ere mutatis mutandis la proc\'edure ainsi mise au point pour obtenir
les deux suites croissantes
$$ \matrix{\displaystyle g_m\subset  g_{m+1}\subset \cdot \cdot \cdot \subset
g_{2m-j}\subset \cdot \cdot \cdot \subset  g_{2m}\, , 1 \leq  j \leq  m ;
\cr\displaystyle h_1\subset  h_2\subset \cdot \cdot \cdot \subset  h_j\subset
\cdot \cdot \cdot \subset  h_m\, . \hfill \cr} $$
Les termes de ces suites sont li\'es par les relations ci-dessous
$$(i) \quad  h_j \subset g_{2m-j} \quad \hbox{pour}\quad 1\leq j \leq m $$
$$(ii) \qquad \qquad[g^o, h_j] \subset h_j.\qquad \qquad \, \, \, \,$$
$$(iii) \qquad ker\omega_{2m-j} = h_j .\qquad \qquad \quad \leqno (9)$$
$$(iv) \qquad [g^o, g_{m-j}] \subset g_{2m-j}.\qquad \qquad$$
$$(v) \quad \omega(h_j \times a, h_j \times a) = 0. $$

On observe en particulier que pour $ j = m $ on a
$$ \omega \left(g_m\times  a,g_m\times a \right) = 0\, . $$
\par
\medskip
\noindent{\bf 2\`eme Etape.} Notons $ G_{2m-j} $ et $ H_j $ les sous-groupes
de Lie connexes d\'etermin\'es
par les sous-alg\`ebres de Lie $ g_{2m-j} $ et $ h_j $ respectivement. Ces
sous-groupes
contenus dans $ G_{2m} $ sont stables par les automorphismes int\'erieurs de $
G $
d\'efinis par les \'el\'ements de $ G^o\, . $ \par
\noindent Cela r\'esulte des relations (8). Cela ach\`eve la preuve du Lemme
1.4.2 . \par
\bigskip
\noindent{\bf Fin de la d\'emonstration du th\'eor\`eme 1.4.1.} \par
A partir des relations (8) on d\'eduit les relations d'inclusions
\'evidentes suivantes
$$ G^o_{n-m+\ell} \subset  H_j\times  G^o_{n-m+\ell} \subset  G_{2m-j}\times
G^o_{n-m+\ell} \subset  G_{2m}\times  G^o_{n-m+\ell} \, . \leqno (10) $$
D'un autre c\^ot\'e la filtration originelle $ F $ vue sous la forme des
produits
semi-directs h\'erit\'es de $ G = G_{2m}\times  G^o $ prend la forme qui suit
$$ G_{2m}\subset  G_{2m}\times  G^o_1 \subset \cdot \cdot \cdot G_{2m}\times
G^o_{n-m+\ell}  \subset \cdot \cdot \cdot G_{2m}\times  G^o\, . \leqno (11) $$
Dans (10) $ G^o_{\ell} $ n'est pas autre chose que le sous-groupe $ H $ dont
l'alg\`ebre de
Lie engendre le noyau de $ \omega \, . $ \par
\noindent On d\'eduit des relations (10) et (11) la filtration $ F_o(G) $
suivante
$$ \doublelow{ G^o_1\subset \cdots G^o_{n-m+\ell} \subset
H_1\times  G^o_{n-m+\ell} \subset \cdots H_m\times  G^o_{n-m+\ell}
\cdots \subset \cr G_{2m}\times  G^o_{n-m+\ell} \subset
G_{2m}\times  G^o_{n-m+\ell +1}\cdots \subset  G_{2m}\times
G^o_{n-m+\ell +k}\subset \cdot \cdot \cdot G_{2m}\times  G^o\, . \hfill \cr}
\leqno (12) $$
Pour voir que le diagramme pond\'er\'e associ\'e \`a (11) est simple, il faut
s'assurer que le graphe orient\'e qui lui est associ\'e ne poss\`ede qu'un
seul
sommet singulier, qui est n\'ecessairement attractif. \par
\noindent En fait on a les diagrammes commutatifs suivants
$$ \matrix{\displaystyle H_{j+1}\times  G^o_{n-m+\ell} & \displaystyle
\leftarrow & \displaystyle H_j\times  G^o_{n-m+\ell} \cr\displaystyle
\downarrow & \displaystyle  & \displaystyle \downarrow \cr\displaystyle
G_{2m-j-1}\times  G^o_{n-m+\ell} & \displaystyle \hookrightarrow &
\displaystyle G_{2m-j}\times  G^o_{n-m+\ell} \cr} \leqno (13) $$
pour $ 0 \leq  j \leq  m $ avec $ H_m= G_m\, ; $ et
$$ \matrix{\displaystyle G^o_{n-m+\ell -k+1} & \displaystyle \leftarrow &
\displaystyle G^o_{n-m+\ell -k} \cr\displaystyle \downarrow & \displaystyle
& \displaystyle \downarrow \cr\displaystyle G_{2m}\times  G^o_{n-m+\ell +k-1}
& \displaystyle \hookrightarrow & \displaystyle G_{2m}\times  G^o_{n-m+\ell
+k} \cr} \leqno (14) $$
pour $ 1 \leq  k \leq  n - m\, . $ Dans les diagrammes (13) et (14) les
fl\`eches sont des
homomorphismes inclusions. Le graphe orient\'e associ\'e \`a (12) a donc un
seul
sommet singulier qui est $ H_m\times  G_{n-m+\ell} \, ; $ voici la portion du
graphe $ d \left(\omega ,F_0 \right) $
qui porte ce sommet singulier
$$ \matrix{\displaystyle \rightarrow & \displaystyle H_{m-1}\times
G^o_{n-m+\ell} & \displaystyle \rightarrow & \displaystyle H_m\times
G^o_{n-m+\ell} & \displaystyle \leftarrow & \displaystyle H_{m-1}\times
G^o_{n-m+\ell} & \displaystyle \leftarrow \cr\displaystyle  & \displaystyle
\downarrow & \displaystyle  & \displaystyle \downarrow & \displaystyle  &
\displaystyle \downarrow & \displaystyle \cr\displaystyle \rightarrow &
\displaystyle H_{m-1}\times  G^o_{n-m+\ell} & \displaystyle \rightarrow &
\displaystyle G_m\times  G^o_{n-m+\ell} & \displaystyle \rightarrow &
\displaystyle G_{m+1}\times  G^o_{n-m+\ell} & \displaystyle \rightarrow \cr} \
. $$
Le diagramme (14) montre que $d (\omega, F_o)$ \`a la forme contract\'e suivante

$ \matrix{\displaystyle O & \displaystyle \rightarrow & \displaystyle O & \displaystyle^{\leftarrow}_{\rightarrow} & \displaystyle O\cr \displaystyle & \displaystyle (n) & \displaystyle & \displaystyle (m) & \displaystyle \cr} $ ;

ainsi le th\'eor\`eme 1.4.1 est d\'emontr\'e . \par
\smallskip
Avant d'entreprendre l'examen des significations g\'eom\'etriques des
propri\'et\'es combinatoires associ\'ees aux graphes $ gr(\omega ,F)\, , $ on
va rappeler un
th\'eor\`eme de $\lbrack$BG 2$\rbrack$ qui donnera au th\'eor\`eme 1.4.1 son int\'er\^et g\'eom\'etrique. \par
\bigskip
\noindent{\bf Th\'eor\`eme }${\bf \lbrack}${\bf BG 2}${\bf \rbrack}$ {\sl Soit
\/}$ G $ {\sl un groupe de Lie compl\`etement r\'esoluble. Soit \/}$ \Gamma $
{\sl un
r\'eseau cocompact dans \/}$ G\, . $ {\sl On suppose que \/}$ \Gamma
\backslash G $ {\sl poss\`ede une structure
kahlerienne \/}$ (\Gamma \backslash G,\omega ,J)\, , $ {\sl alors la classe de
Kahler \/}$ [\omega] \in  H^2_{Dk}(\Gamma \backslash G) $ {\sl contient
une forme symplectique \/}$ G ${\sl -homog\`ene \/}$ \omega_ o $ {\sl avec les
propri\'et\'es suivantes\/} \par
\medskip
\noindent{\sl (i) \/}$ \omega_ 0 $ {\sl est la projection d'une forme
invariante \`a gauche not\'ee encore
\/}$ \omega_ 0\in  \Omega^ 2_{\ell}( G), $ $\lbrack$c'est la $ G
$-homog\'en\'eit\'e$\rbrack$. \par
\medskip
\noindent{\sl (ii) La restriction \`a l'id\'eal d\'eriv\'e \/}$ [g,g] $ {\sl
de \/}$ \omega_ 0 $ {\sl est une 2-forme
symplectique\/} . \par
\smallskip
\noindent Le th\'eor\`eme 2 de $\lbrack$BG 2$\rbrack$ est plus complet que
celui ci-dessus. On d\'eduit en
fait de (ii) que l'orthogonal symplectique de $ [g,g] $ est une sous alg\`ebre
de Lie ab\'elienne $ a $ de $ g\, . $ Il r\'esulte des propositions 1.3.1 et
1.3.2 que $ g $
poss\`ede une filtration dont le diagramme associ\'e est semi-simple et
semi-nilpotent. \par
\noindent Le th\'eor\`eome 1.4.1 assure que $ \omega_ 0 $ poss\`ede un
diagramme pond\'er\'e simple
$$ {\rm O} \rightarrow {\rm O}^{\leftarrow}_{ \rightarrow} {\rm O} \, . $$
Soit $ G $ un groupe de Lie compl\`etement r\'esoluble contenant un r\'eseau
cocompact $ \Gamma \, . $ Soit $ \omega  \in  \Omega^ 2_{\ell}( G) $ une forme
symplectique. On note $ (\Gamma \backslash G,\omega) $ la
vari\'et\'e symplectique obtenue par la projection canonique de $ G $ sur $
\Gamma \backslash G\, . $
Gr\^ace au th\'eor\`eme $\lbrack$BG 2$\rbrack$ ci-dessus on d\'eduit du
th\'eor\`eme 1.4.1 le \par
\bigskip
\noindent{\bf Corollaire 1.4.2.} {\sl Si la vari\'et\'e \/}$ \Gamma
\backslash G $ {\sl poss\`ede une structure kahlerienne dont
la forme de Kahler \/}$ \omega \, $ est $G$-homog\`ene,  {\sl alors cette derni\`ere
poss\`ede un diagramme
pond\'er\'e simple \/}$ {\rm O} \rightarrow {\rm O}^{\leftarrow}_{
\rightarrow} {\rm O} \ . $ \par
\smallskip
\noindent Naturellement il revient au m\^eme de dire que $ \omega $ poss\`ede
un graphe orient\'e
$ gr(\omega ,F) $ ayant un seul sommet singulier . \par
\skiplines 2
\noindent{\bf 2. DYNAMIQUES DIFF\'ERENTIABLES ET DIAGRAMMES POND\'ER\'ES}.
\par
Ce \S 2 est consacr\'e \`a deux questions. L'une est de nature
topologique ; l'autre question concerne la r\'eductibilit\'e des \lq\lq paires
transitives\rq\rq\ au sens de A.L. ONISHCHIK, $\lbrack$0N$\rbrack$. Dans la
suite on entend par
syst\`emes dynamiques diff\'erentiables les actions diff\'erentiables des
groupes de Lie. Le probl\`eme topologique qui nous int\'eresse est celui de
l'existence de feuilletages lagrangiens invariants par des syst\`emes
dynamiques symplectiques. En situation compl\`etement r\'esoluble ce
pro\-bl\`eme
est en rapport avec la nature des diagrammes pond\'er\'es $ d(\omega ,F) $ ou
si l'on
veut des graphes orient\'es associ\'es, $\lbrack$au sens des d\'efinitions
1.3.3 et
1.3.4). Le probl\`eme de \lq\lq r\'eductibilit\'e\rq\rq\ des paires
transitives est quant \`a
lui li\'e \`a la nature des sommets singuliers des diagrammes pond\'er\'es.
Nous
allons pr\'eciser ces liaisons une fois rappel\'ees les notions de degr\'e primitif et de
sym\'etrie des vari\'et\'es homog\`enes, voir $\lbrack$HW$\rbrack$ pour le
cas des groupes
compactes, notions par ailleurs \'etroitement li\'ees aux poids des sommets
singuliers des diagrammes pond\'er\'es. \par
\bigskip
\noindent{\bf 2.1. Degr\'es primitifs des vari\'et\'es homog\`enes.} \par
Les vari\'et\'es consid\'er\'ees sont diff\'erentiables de classe $
C^{\infty} \, , $
connexes et d\'enombrables \`a l'infini. Une vari\'et\'e $ M $ est dite
homog\`ene quand
elle poss\`ede un {\sl groupe de Lie des transformations transitif.\/} La
collection
constitu\'ee des groupes de Lie des transformations transitifs de $ M $ sera
d\'esign\'ee par $ {\cal L}(M)\, . $ Tout $ G \in  {\cal L}(M) $ permet
l'identification de $ M $ avec une
paire $ (H,G) $ o\`u $ H $ est le sous-groupe de $ G $ stabilisateur d'un
point dans $ M\, . $
Posons
$$ d^G(M) = \frac{dim H}{codim H +1}$$
$$ d(M) = min \{ d^G(M) \}. $$
$$ G \in {\cal L}(M)$$
\par
\smallskip
\noindent{\bf D\'efinition 2.1.1.}  {\sl Un couple $ (H,G)$ est dit quasi-primitif ( respectivement primitif) si $G$ ne contient aucun sous-groupe (repectivement aucun sous-groupe distingu\'e) $ G_o$  avec $G_o \in \cal {L}(G/H)$
\/}  \par
\bigskip
\noindent Classiquement $[H W]$ on appelle degr\'e de sym\'etrie de $ M $ le nombre
entier
$$ \doublelow{ \delta( M) = \sup(\dim(G))\, , \cr G \in  {\cal L}(M)
\quad et \quad G \quad est \quad compact .\cr} $$
Nous allons maintenant nous restreindre \`a la sous-collection $ {\cal R}(M)
\subset  {\cal L}(M) $
constitu\'ee des groupes de Lie r\'esolubles. Autrement dit $ M $ est une
solvari\'et\'e (solvmanifolds) si $ {\cal R}(M) $ n'est pas
vide. Pour
une \'etude d\'etaill\'ee des $ {\cal R}(M) $ lorsque $ M $ est compacte, on
renvoie au
panorama de L. Auslander $ \left[AL \right] $ et \`a $ \left[MO \right] $
et $ \left[MO \right] $ pour des consid\'erations
plus cibl\'ees, $\lbrack$e.g. cohomologie ou Homotopie des dites
vari\'et\'es$\rbrack$. \par
\noindent Les vari\'et\'es ayant le degr\'e primitif nul sont les
quotients des groupes
de Lie par leurs sous-groupes ferm\'es discrets. Quand ces vari\'et\'es sont
compactes, les sous-groupes de Lie stabilisateurs des points de vari\'et\'es en question sont des r\'eseaux
cocompactes, $\lbrack$RA$\rbrack$. \par
\noindent La notion de degr\'e qui est \'etroitement li\'ee aux notions
combinatoires des
num\'eros 1 et 2 est celle de degr\'e relatif. \par
\noindent On fixe un groupe de Lie $ G \in  {\cal L}(M)\, , $ ou ce qui
revient au m\^eme une
g\'eom\'etrie particuli\`ere dans la vari\'et\'e $ M\, . $ On d\'esigne par $
{\cal L}_G(M) $ la
sous-collection des $ G^\prime \in  {\cal L}(M) $ avec $ G^\prime \subset
G\, . $ La question est d'\lq\lq optimiser\rq\rq\
dans un sens \`a pr\'eciser le choix de $ G^\prime $ dans $ {\cal L}_G(M)\, .
$ Puisque pour chaque
choix $ G^\prime \in  {\cal L}_G(M)\, , $ $ M $ est identifi\'e \`a une paire
$ (H^\prime ,G^\prime) \, , $ les choix des $ H^\prime $
\'etant param\'etr\'es par les points $ x $ de $ M\, ; $ on associera \`a la
paire $ (H^\prime ,G^\prime) $
le nombre rationnel
$$ r_{G^{\prime}}(M) = \frac{dim H^{\prime}}{codim H^{\prime} +1} = \frac{dim H^{\prime}}{dim M +1} ; $$

o\`u $M$ est l'espace homog\`ene $ G^\prime/H^\prime $ .
Ce nombre ne d\'epend que de $ G^\prime \in  {\cal L}_G(M)\, ; $ de
sorte que pour $ G^\prime $ et $ G^{\prime\prime} $ dans $ {\cal L}_G(M)\, , $
la comparaison des paires $ (H^\prime ,G^\prime) $ et $ (H^{\prime\prime}
,G^{\prime\prime}) $ peut en fait d\'ependre des
invariants topologiques ou g\'eom\'etriques du couple $ (G,M)\, , $ (voir par
$\lbrack$ON$\rbrack$
pour le cas des vari\'et\'es riemanniennes homog\`enes compactes). Notre but est de
comparer $ (H^\prime ,G^\prime) $ et $ (H^{\prime\prime} ,G^{\prime\prime}) $
au moyen des nombres rationnels $ r_{G^\prime}( M) $ et
$ r_{G^{\prime\prime}}( M)\, . $ Posons
$$ d_G(M) = \doublelow{ \min \{\frac{dim H^{\prime}}{codim H^{\prime}} \} \cr G^\prime \in {\cal L}_G(M) \cr} $$
\par
\smallskip
\noindent{\bf D\'efinition 2.1.2.} {\sl Le nombre rationnel \/}$ d_G(M) $
{\sl est appel\'e le \/}degr\'e primitif de  $ (G,M) . $ \par
\bigskip
\noindent Il est clair que si $ G^\prime $ et $ G^{\prime\prime} $ sont dans $
{\cal L}_G(M) $ la relation d'inclusion $ G^{\prime\prime}  \subseteq
G^\prime $
entraine \`a $ d_{G^{\prime\prime}}( M) \leq  d_{G^\prime}( M)\, . $ \par
Exemples. a) Si une vari\'et\'e $M$ poss\`ede une structure de groupe de lie compatible avec sa topologie g\'en\'erale alors $d(M) = 0$.\\
 b) Si $\Gamma$ est un sous-groupe discret ferm\'e dans un groupe de lie connexe $G$ alors $d({\Gamma}\backslash G) = 0 .$
\noindent  On a l'in\'egalit\'e $ d(M) \leq  d_G(M) $ pour tout $
G \in  {\cal L}(M)\, . $ \par
\bigskip
\noindent{\bf D\'efintion 2.1.3.} {\sl Pour \/}$ G \in  {\cal L}(M) $ {\sl on
dit que \/}$ (G,M) $ {\sl est  quasi-primitifs si
 \/}$ {\cal L}_G(M) = \{ G\} \, ; $ {\sl un groupe de lie \/}$ G \in  {\cal
L}(M) $ {\sl est dit non quasi-primitif si le cardinal
de \/}$ {\cal L}_G(M) $ {\sl est sup\'erieur \`a 1.\/} \par
\bigskip
\noindent On voit sans difficult\'e que les propri\'et\'es suivantes sont
\'equivalentes \par
\medskip
\noindent (i) $ G \in  {\cal L}(M) $ {\sl est  quasi-primitif ;\/} \par
\medskip
\noindent (ii) $ d_G(M) = r_G(M)\  $ \par
\bigskip
{\bf 2.2. R\'eductibilit\'e des diagrammes
pond\'er\'es.} \\
D\'esormais on se restreind \`a la cat\'egorie des groupes de Lie
compl\`etement r\'esolubles et \`a la dynamique de ces groupes dans leurs
vari\'et\'es symplectiques homog\`enes. \par
Soit $ G $ un groupe de Lie compl\`etement r\'esoluble. Pour une forme
 $ \omega  \in  \Omega^ 2_{\ell}( G)\, , $ on consid\`ere la
vari\'et\'e symplectique $ (G/H,\omega) $ o\`u $ H $ est
le sous-groupe connexe dont l'alg\`ebre de Lie engendre le noyau de $ \omega
\, . $ On
note encore $ \omega $ la forme symplectique de $ G/H $ h\'erit\'ee de $
\omega  \in  \Omega^ 2_{\ell}( G)\, . $
Naturellement $ G \in  {\cal L}(M) $ o\`u $ M $ d\'esigne la vari\'et\'e $
G/H\, . $ On identifie $ M $ avec
la paire $ (H,G) $ et le couple $ (M,\omega) $ avec le triplet $ (H,G,\omega)
$ ou $ (G,H,\omega) \, . $ \par
\smallskip
\noindent Les notations $ {\cal L}(M,\omega) $ et $ {\cal L}_G(M,\omega) $ ont
un sens \'evident. On peut parler des
degr\'es primitifs symplectique et de degr\'e primitif relatif.
Nous les d\'esignerons respectivement par
$$\doublelow{d(M, \omega) = min (d^{G^{\prime}}(M)) , \cr G^{\prime} \in {\cal L}_{G}(M, \omega)\cr}$$
$$\doublelow{d_{G}(M, \omega) = min (r_{G^{\prime}}(M)) , \cr G^{\prime} \in {\cal L}(M, \omega)\cr}$$
Soit maintenent $ (\omega , F) \in \Omega^ 2_{\ell}( G) \times {\cal F}(G) $ et soit $ \left(G_k,H_k,\omega_ k
\right) $ le syst\`eme associ\'e. On sait
pouvoir identifier $ M_k $ avec $ \left(G_k,H_k \right) $ et $
\left(M_k,\omega_ k \right) $ avec $ \left(G_k,H_k,\omega_ k \right)\, . $ Par
abus
de notations habituel $ \omega_ k $ d\'esigne tout aussi bien $ \omega_ k\in
\Omega^ 2_{\ell} \left(G_k \right) $ que la forme
symplectique de $ M_k $ qui en est d\'eduite. \par
\noindent Nous allons \'etablir des liens entre les natures des sommets du
graphe
orient\'e $ gr(\omega ,F)\, , $ (ou du diagramme pond\'er\'e $ d(\omega ,F)) $
avec la primitivi\'e
de $ G_k $ dans $ {\cal L}_{G_k} \left(M_k,\omega_ k \right) = {\cal L}_{G_k}
\left(M_k \right)\, . $ \par
\noindent Les sommets r\'eguliers sont en fait de deux types : Ceux de type $
(\Irm) $ sont
appel\'es des sommets hyperboliques croissants, alors que ceux de type $
(\IIrm) $
sont de type hyperbolique d\'ecroissante $\lbrack$NB 5$\rbrack$. De sorte que
dans les
coordonn\'ees $ \left(G_k,H_k \right) $ avec les relations d'inclusion mises
en \'evidence au
num\'ero 1 les sommets r\'eguliers et singuliers sont pr\'esent\'es de fa\c
con
dynamique comme ci-dessous, (voir $\lbrack$NB 5$\rbrack$).
$$ \matrix{\includegraphics{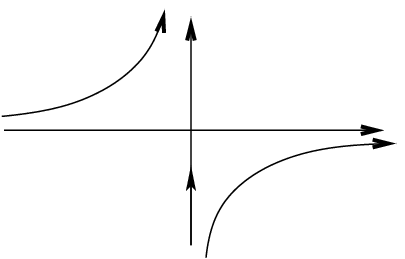}\displaystyle  & \displaystyle {\rm pour} & \displaystyle {\rm
O}^{\rightarrow}_{ \rightarrow} \doubleup{ S \cr {\rm O} \cr}^{
\rightarrow}_{ \rightarrow} {\rm O} \cr\displaystyle  & \displaystyle  &
\displaystyle \cr\includegraphics{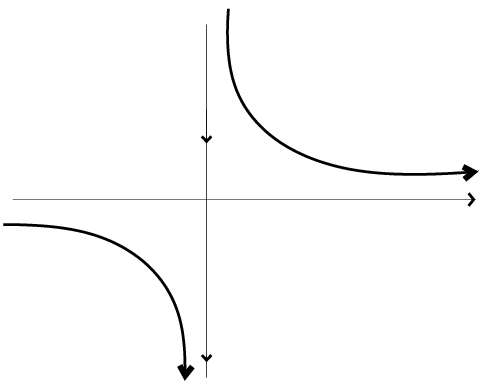}\displaystyle  & \displaystyle {\rm pour} & \displaystyle
{\rm O}^{\leftarrow}_{ \rightarrow} \doubleup{ S \cr {\rm O} \cr}^{
\leftarrow}_{ \rightarrow} {\rm O} \cr\displaystyle  & \displaystyle  &
\displaystyle \cr\includegraphics{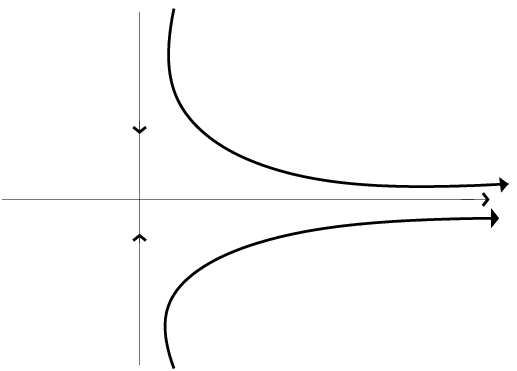}\displaystyle  & \displaystyle {\rm pour} & \displaystyle
{\rm O}^{\rightarrow}_{ \rightarrow} \doubleup{ S \cr {\rm O} \cr}^{
\leftarrow}_{ \rightarrow} {\rm O} \cr\displaystyle  & \displaystyle  &
\displaystyle \cr\includegraphics{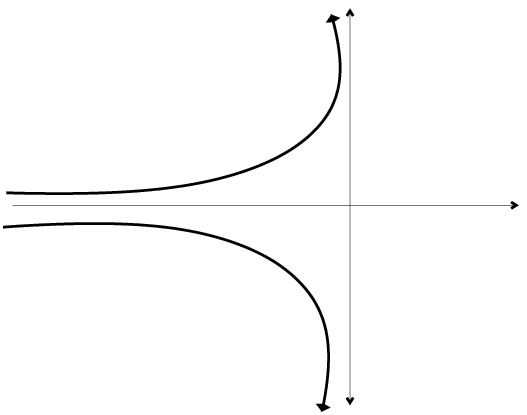}\displaystyle  & \displaystyle {\rm pour} & \displaystyle
{\rm O}^{\leftarrow}_{ \rightarrow} \doubleup{ S \cr {\rm O} \cr}^{
\rightarrow}_{ \rightarrow} {\rm O\, .} \cr} \leqno (15) $$
\noindent \overfullrule=0mm \par
\noindent La d\'efinition suivante est l'analogue des relations de ONISHCHIK
pour les
paires riemaniennes transitives $\lbrack$ON$\rbrack$. \par
\noindent Deux sommets $ S_k $ et $ S_{\ell} $ d'un graphe $ {\rm gr}(\omega
,F) $ sont dites {\sl\'equivalentes\/} si les
vari\'et\'es $ G_k/H_k $ et $ G_{\ell} /H_{\ell} $ sont canoniquement identiques, ils
sont{\sl\ in\'equivalents\/} sinon. L'expression canoniquement identiques a le sens suivant : $ dim G_k / H_k = dim G_{\ell}/H_{\ell}$ et soit $G_k \in {\cal L}(G_{\ell}/H_{\ell})$ soit $ G_{\ell} \in {\cal L}(G_k /H_k)$ . \par
\noindent Lorsque deux sommets $ S_k $ et $ S_{\ell} $ sont \'equivalents on a
n\'ecessairement l'une
des relations d'inclusion entre paires $ \left(G_k,H_k \right) \subset
\left(G_{\ell} ,H_{\ell} \right) $ ou
$ \left(G_{\ell} ,H_{\ell} \right) \subset  \left(G_k,H_k \right) , $ (voir
$\lbrack$O$\rbrack$). Lorsque $ \left(G_k,H_k \right) $ est \'equivalent \`a $
\left(G_{\ell} ,H_{\ell} \right) $
avec $ \left(G_k,H_k \right) \subset  \left(G_{\ell} ,H_{\ell} \right)\, , $
la paire $ \left(G_{\ell} ,H_{\ell} \right) $ est {\sl r\'eductible \`a la
paire\/} $ \left(G_k,H_{\ell} \right)\, ; $
dans ce cas tous les sommets situ\'es entre $ \left(G_k,H_k \right) $ et $
\left(G_{\ell} ,H_{\ell} \right) $ sont
r\'eguliers de type hyperbolique croissant. \par
\noindent Les affirmations qui suivent sont faciles \`a v\'erifier. \par
\bigskip
\noindent{\bf Proposition 2.2.1.} \par
\noindent (i) {\sl Un sommet hyperbolique croissant est \'equivalent aux deux
sommets qui
lui sont adjacents.\/} \par
\noindent{\sl (ii) Un sommet hyperbolique d\'ecroissant est non \'equivalent
aux deux
sommets qui lui sont adjacents.\/} \par
\noindent{\sl (iii) Un sommet singulier attractif est \'equivalent au sommet
qui lui est
adjacent \`a gauche.\/} \par
\noindent{\sl (iv) Un sommet singulier r\'epulsif est \'equivalent au sommet
qui lui est
adjacent \`a droite }. \par
\bigskip
\noindent Nous nous proposons d'\'etablir un lien \'etroit entre le degr\'e
primitif relatif et la combinatoire des graphes $\lbrack$ou si l'on veut, des
diagrammes
pond\'er\'es). \par
Soit $ G $ un groupe de Lie compl\`etement r\'esoluble. Soit $ \omega  \in
\Omega^ 2_{\ell}( G) $ . Soit $ H $ le sous-groupe connexe de $ G $ dont l'alg\`ebre de
Lie
engendre le noyau de $ \omega \, . $ La vari\'et\'e symplectique $
(G/H,\omega) $ est identifi\'ee au
triplet $ (G,H,\omega) \, . $ \par
\bigskip
\noindent{\bf Th\'eor\`eme 2.2.2.} {\sl On suppose que la paire \/}$ (G,H) $
{\sl est quasi-primitive \/}$ r_G(M) = d_G(M))\, . $ {\sl Si un diagramme pond\'er\'e \/}$ d(\omega ,F)
$ {\sl est connexe, alors il
poss\`ede au plus trois sommets singuliers }. \par
\bigskip
\noindent{\sl Esquisse de d\'emonstration.\/} \par
La preuve s'appuie sur les propositions 1.3.1 et 1.3.2. En fait
l'hypoth\`ese $ d_G(M) = r_G(M) $ entra\^\i ne que le sommet $ (G,H) = S $ de
$ d(\omega ,F) $ est
hyperbolique d\'ecroissante \`a sa gauche. Puisque $ d(\omega ,F) $ est
connexe il a la
forme
$$ {\rm O} \rightarrow {\rm O} \cdot \cdot \cdot {\rm O}^{\leftarrow}_{
\rightarrow} \doubleup{ S \cr {\rm O\, .} \cr} $$
Soit $ S^\prime $ le sommet singulier le plus proche du sommet $ S\, ; $ ce
sommet est
attractif. Si $ d(\omega ,F) $ est simple alors $ S^\prime $ est l'unique
sommet singulier de
$ d(\omega ,F)\, . $ on d\'eduit que $ d(\omega ,F)$ poss\`ede un seul sommet singulier. \par
\noindent Si ce n'est pas le cas la vari\'et\'e symplectique $ M^\prime $
d\'efinie par le sommet
singulier $ S^\prime $ est de dimension $ > 0\, . $ La proposition 1.3.2
assure que le
sommet $ S^\prime $ est \'equivalent \`a un sommet singulier $
S^{\prime\prime} $ qui est nilpotent et
r\'epulsif. La connexit\'e de $ d(\omega ,F) $ joints \`a la proposition 1.3.1
entra\^\i ne
qu'il existe un et un seul sommet singulier $ S^{\prime\prime\prime} $ \`a
gauche de $ S^{\prime\prime} \, . $ Cela
ach\`eve la d\'emonstration du th\'eor\`eme . \par
\medskip
\noindent{\sl Remarque 2.2.3.\/} Du theor\`eme 2.2.2 on d\'eduit qu'en fait
un diagramme
connexe poss\`ede au plus quatre sommets singuliers. Avec la convention de
contraction sommets hyperboliques voici les seuls mod\`eles des diagrammes
connexes

$$(\alpha)\quad \matrix { \displaystyle O \displaystyle\rightarrow\displaystyle O \displaystyle ^{\leftarrow}_{\rightarrow} \displaystyle O \displaystyle \rightarrow \displaystyle O \displaystyle ^{\leftarrow}_{\rightarrow} \displaystyle O \displaystyle \rightarrow \displaystyle O\cr\displaystyle \displaystyle (m_1)\displaystyle \displaystyle (n_1) \displaystyle \displaystyle (m_2) \displaystyle \displaystyle (n_2) \displaystyle \displaystyle (m_3) \displaystyle\cr } {\rm associ\acute e\ au\ graphe} \matrix {\displaystyle O&\displaystyle\rightarrow &\displaystyle O&\displaystyle \leftarrow &\displaystyle O&\displaystyle \rightarrow &\displaystyle O&\displaystyle\leftarrow &\displaystyle O&\displaystyle\rightarrow &\displaystyle O\cr\displaystyle\searrow &\displaystyle &\displaystyle\swarrow &\displaystyle &\displaystyle\searrow &\displaystyle &\displaystyle\swarrow &\displaystyle &\displaystyle \searrow &\displaystyle &\displaystyle\swarrow\cr\displaystyle &\displaystyle O &\displaystyle &\displaystyle &\displaystyle &\displaystyle O &\displaystyle &\displaystyle &\displaystyle &\displaystyle O &\displaystyle\cr } $$

$$(\beta)\quad \matrix { \displaystyle O \displaystyle\rightarrow\displaystyle O \displaystyle ^{\leftarrow}_{\rightarrow} \displaystyle O \displaystyle \rightarrow \displaystyle O \displaystyle ^{\leftarrow}_{\rightarrow} \displaystyle O\cr\displaystyle \displaystyle (m_1)\displaystyle \displaystyle (n_1) \displaystyle \displaystyle (m_2) \displaystyle \displaystyle (n_2) \displaystyle \cr } {\rm \quad associ\acute e\ au\ graphe \quad } \matrix {\displaystyle O&\displaystyle\rightarrow &\displaystyle O&\displaystyle \leftarrow &\displaystyle O&\displaystyle \rightarrow &\displaystyle O&\displaystyle\leftarrow &\displaystyle O\cr\displaystyle\searrow &\displaystyle &\displaystyle\swarrow &\displaystyle &\displaystyle\searrow &\displaystyle &\displaystyle\swarrow &\displaystyle &\displaystyle \downarrow\cr\displaystyle &\displaystyle O &\displaystyle &\displaystyle &\displaystyle &\displaystyle O &\displaystyle &\displaystyle &\displaystyle O\cr } $$

$$(\gamma)\qquad \matrix { \displaystyle O \displaystyle\rightarrow\displaystyle O \displaystyle ^{\leftarrow}_{\rightarrow} \displaystyle O \displaystyle \rightarrow \displaystyle O \cr\displaystyle \displaystyle (m_1)\displaystyle \displaystyle (n_1) \displaystyle \displaystyle (m_2) \displaystyle \cr } {\rm \quad associ\acute e\ au\ graphe \quad } \matrix {\displaystyle O&\displaystyle\rightarrow &\displaystyle O&\displaystyle \leftarrow &\displaystyle O&\displaystyle \rightarrow &\displaystyle O\cr\displaystyle\searrow &\displaystyle &\displaystyle\swarrow &\displaystyle &\displaystyle\searrow &\displaystyle &\displaystyle\swarrow\cr\displaystyle &\displaystyle O &\displaystyle &\displaystyle &\displaystyle &\displaystyle O &\displaystyle \cr } $$

$$(\delta)\qquad \matrix { \quad \displaystyle O \displaystyle\rightarrow\displaystyle O \displaystyle ^{\leftarrow}_{\rightarrow} \displaystyle O \cr\displaystyle \displaystyle (m_1)\displaystyle \displaystyle (n_1) \displaystyle \cr } \qquad \qquad {\rm \quad associ\acute e\ au\ graphe \quad } \matrix {\displaystyle O&\displaystyle\rightarrow &\displaystyle O&\displaystyle \leftarrow &\displaystyle O\cr\displaystyle\searrow &\displaystyle &\displaystyle\swarrow &\displaystyle &\displaystyle\downarrow \cr\displaystyle &\displaystyle O &\displaystyle &\displaystyle  &\displaystyle O \quad \cr } \quad $$
Il faut comprendre l'association de diagramme au graphe ci-dessus de la fa\c{c}on qui suit.\\
($\alpha$) le diagramme poss\`ede quatre sommets singuliers; tous les sommets r\'eguliers compris entre deux sommets singuliers cons\'ecutifs sont deux \`a deux \'equivalents, ces sommets r\'eguliers sont r\'epartis en trois classes d'\'equivalence dont chacune contient plusieurs sommets r\'eguliers.\\
($\beta$) Le diagramme poss\`ede deux sommets singuliers ; il y a trois classes d'\'equivalence de sommet r\'egulier; la classe situ\'ee \`a l'extr\^eme droite contient un seul sommet r\'egulier.\\
($\gamma$) Le diagramme poss\`ede deux sommets singuliers et il y a deux classes de sommet r\'egulier dont chacune contient plusieurs sommets r\'eguliers.\\
($\delta$) Le diagramme poss\`ede un seul sommet singulier, il y a deux classes d'\'equivalence de sommet r\'egulier, la classe situ\'ee \`a l'extr\^eme droite contient un seul sommet r\'egulier.\\
Les cas ($ \alpha$) et ($ \gamma$) sont des diagrammes et graphes des paires
non quasi-primitives les cas ($ \beta$) et ($ \gamma$) sont ceux des paires quasi-primitives. De plus ($ \gamma$) est simple \par
\bigskip
\noindent {\bf 2.3. Sommets singuliers et feuilletages
invariants.} \par
Soit $ (M,\omega) $ une vari\'et\'e symplectique homog\`ene et $ {\cal
L}(M,\omega) $ la
collection de ses dynamiques diff\'erentiables. Supposons que $ {\cal
L}(M,\omega) $
contient un groupe de Lie compl\`etement r\'esoluble $ G\, . $ On sait que
chaque
choix d'un point $ x $ dans $ M $ donne lieu \`a l'identification de $ M $
avec une
paire $ (G,H) $ et de $ (M,\omega) $ avec un triplet $ (G,H,\omega) \, . $ A
chaque $ F \in  {\cal F}(G) $
correspond les objets combinatoires $ {\rm gr}(\omega ,F) $ et $ d(\omega ,F)
$ ayant pour {\sl sommet
extr\'emal\/} $ (G,H) = S\, . $ \par
\noindent Nous allons montrer qu'en fait les poids des sommets r\'epulsifs des
$ d(\omega ,F) $
constituent des obstructions \`a l'existence de feuilletages lagrangiens
invariants par $ G\, . $ \par
\bigskip
\bigskip
\noindent{\bf 2.3.0.  Simplicit\'e des graphes }$ {\rm gr}(\omega ,F) $ {\bf et
feuilletages lagrangiens.} \par
Les donn\'ees $ (M,\omega)  , $ $ G \in  {\cal L}(M,\omega) $ sont comme
ci-dessus, viz $ G $ est
compl\`etement r\'esoluble. On identifie $ M $ et $ (M,\omega) $ avec $ (G,H)
$ et $ (G,H,\omega) $
respectivement. Nous conservons l'abus d'\'ecriture convenu, viz $ \omega $
est
l'image inverse de la forme symplectique de $ (M,\omega) $ par
l'identification
$ G/H \sim  M\, . $ \par
On sait construire dans $ M $ d'autres g\'eom\'etries li\'ees \`a $
(M,\omega) \, . $ E.g : \par
\noindent (i) Feuilletages lagrangiens, $\lbrack$NB 1$\rbrack$. \par
\noindent (ii) Structures bilagrangiennes, $\lbrack$NB 5$\rbrack$. \par
\noindent (iii) Structures k\"ahleriennes affinement plates $\lbrack$NB
4$\rbrack$. \par
En g\'en\'eral ces structures ne sont pas des invariants de la
dynamique donn\'ee $ G \in  {\cal L}(M,\omega) \, , $ bien que celle-ci joue
un grand r\^ole dans
leur construction, $\lbrack$NB 5$\rbrack$. A l'oppos\'e de la derni\`ere
remarque, il arrive
que l'existence seule d'un certain type de dynamique $ G \in  {\cal L}(M) $
ait des
cons\'equences spectaculaires sur la g\'eom\'etrie globale de la vari\'et\'e $
M\, . $ Les
th\'eor\`emes de  $\lbrack$BG 1$\rbrack$ et de Dusa Mc Duff, $\lbrack$Mc
D$\rbrack$ en sont des exemples
frappants. La conjecture des tores plats est l'un des motivations de ce
travail . \par
\noindent Ci-dessous un premier r\'esultat qui lie la simplicit\'e des $ {\rm
gr}(\omega ,F) $ aux
feuilletages lagrangiens invariants. On va commencer par fixer quelques
identifications combinatoires. \par
Soit $ G $ un groupe de Lie compl\`etement r\'esoluble. Fixons une 2-forme
 $ \omega  \in  \Omega^ 2_{\ell}( G)\, . $ Consid\'erons deux s\'eries
de compositions $ F $ et $ F^\prime $ dans
$ {\cal F}(G)\, . $ Nous allons choisir des \'el\'ements de comparaison de $
{\rm gr}(\omega ,F) $ \`a $ {\rm gr}(\omega ,F^\prime) $
et de $ d(\omega ,F) $ \`a $ d(\omega ,F^\prime) $ respectivement. Si on fait
abstraction de
sous-groupes composant $ F $ et $ F^\prime $ les graphes (resp. les
diagrammes) associ\'es
\`a $ (\omega ,F) $ et \`a $ (\omega ,F^\prime) $ peuvent pr\'esenter des
configurations g\'eom\'etriques ou
combinatoires identiques. Nous allons nous int\'eresser aux sommets
singuliers. \par
\bigskip
\noindent{\bf D\'efinition 2.3.1.} {\sl Deux graphes \/}$ {\rm gr}(\omega ,F)
$ {\sl et \/}$ {\rm gr}(\omega ,F^\prime) $ $\lbrack${\sl resp. deux
diagrammes
pond\'er\'es \/}$ d(\omega ,F) $ {\sl et \/}$ d(\omega ,F^\prime) ) $ {\sl
sont dits \'equivalents si et seulement s'ils
ont des sommets singuliers identiques.\/} \par
\bigskip
\noindent{\sl Remarque 2.3.2.\/} \par
\noindent Avec la convention de contraction des sommets r\'eguliers adapt\'ee
deux graphes distincts peuvent donner lieu au m\^eme diagramme contract\'e.
 Ainsi {\sl aller des classes d'isomorphisme des
alg\`ebres de
Lie semi-simples sur un corps alg\'ebriquement clos\/} \`a l'unique {\sl
diagramme\/} de
Dynkin associ\'e est l'analogue de passage des classes d'\'equivalence des
graphes $ {\rm gr}(\omega ,F)\, , $ $ F \in  {\cal F}(G) $ aux diagrammes
contract\'es . \par
\smallskip
\noindent On notera que les
multiplicit\'es des arr\^etes ( ou sommets) hyperboliques des diagrammes contract\'es permet de distinguer les graphes simples;  par exemple le diagramme
$$ {\rm O} \doublelow{ \rightarrow \cr (1) \cr} \doubleup{ S \cr {\rm O} \cr}
\doublelow{^{ \leftarrow}_{ \rightarrow} \cr (1) \cr} {\rm O} $$
correspond \`a une surface orient\'ee; son graphe est distinct de celui du diagramme
$$ {\rm O} \doublelow{ \rightarrow \cr (2) \cr} {\rm O} \doublelow{^{
\leftarrow}_{ \rightarrow} \cr (2) \cr} {\rm O} $$
qui correspond \`a une vari\'et\'e symplectique de dimension 4. Par
cons\'equent
l'expression \lq\lq diagramme contract\'e\rq\rq\ s'entendra toujours
multiplicit\'es des
sommets r\'eguliers comprises m\^eme si on omet volontairement de les
expliciter. \par
\noindent Les notations $ G,\omega ,(M,\omega)  = (G,H,\omega) $ sont celles
fix\'ees dans les sections et sous-sections pr\'ec\'edentes. \par
\bigskip
\noindent{\bf Th\'eor\`eme 2.3.3.} {\sl Il y a correspondance bijective entre
les classes
d'\'equivalence des graphes orient\'es simples \/}$ {\rm gr}(\omega ,F) $
{\sl dans \/}$ \{ {\rm gr}(\omega ,F)/F\in{\cal F}(G)\} $
{\sl et les feuilletages lagrangiens dans \/}$ (M,\omega) $ {\sl qui sont
invariants par \/}$ G\  $ \par
\bigskip
\noindent{\sl D\'emonstration.\/} Soit $ {\rm gr}(\omega ,F) $ un graphe
orient\'e simple. Si $ 2m $ est la
dimension de $ G/H $ le diagramme $ d(\omega ,F) $ poss\`ede exactement $ m $
sommet
hyperboliques d\'ecroissants et un seul sommet singulier attractif $ S\, . $
Il a
la forme contract\'ee
$$ {\rm O} \rightarrow {\rm O}^{\leftarrow}_{ \rightarrow} {\rm O} \ .  $$
Le sommet singulier $ S $ a pour coordonn\'ees la paire $ S =
\left(G_{m+\ell} ,H \right) $ o\`u
$ \dim H = \ell \, . $ Bien entendu on a $ H \subset  G_{m+\ell} $ et $
i^\ast_{ m+\ell} \omega  = 0\, . $ Le feuilletage de $ G $
par les classes \`a droite $ \gamma G_{m+\ell} $ est invariant par les
translations \`a gauche
dans $ G $ et est projet\'e en un feuilletage lagrangien dans $ (M,\omega) \,
. $ Le
feuilletage ainsi obtenu est invariant par l'action de $ G $ dans $
(M,\omega) \, ; $ il
ne d\'epend d'autre part que de la classe d'\'equivalence de $ {\rm
gr}(\omega ,F)\  $ \par
\smallskip
\noindent R\'eciproquement, supposons que l'on ait dans $ (M,\omega) $ un
feuilletage
lagrangien $ L $ qui est invariant par $ G\, . $ Soit $ x_0\in  M $ un point
dont le
sous-groupe stabilisateur est $ H\, . $ Notons $ \pi $ la projection orbitale
$ \gamma  \mapsto  \gamma x_0 $
et $ L^\prime  = \pi^{ -1}(L)\, ; $ ce dernier est un feuiletage invariant par
les
translations \`a gauche dans $ G\, . $ La feuille $ L^\prime( e) $ de $
L^\prime $ qui contient
l'\'el\'ement neutre est un sous-groupe de $ G $ de dimension $ m + \ell $ et
contenant $ H\, . $
Puisque $ G $ est compl\`etement r\'esoluble, il r\'esulte du classique
deuxi\`eme
th\'eor\`eme de Lie que l'on a un $ F \in  {\cal F}(G) $ tel que la feuille $
L^\prime( e) $ soit un
sommet du graphe $ {\rm gr}(\omega ,F)\, . $ \par
\noindent Puisque $ H $ est inclus dans $ L^\prime( e) $ et que $ i^\ast_{
L^\prime( e)}\omega  = 0\, , $ $ L^\prime( e) $ est l'unique
sommet singulier de $ d(\omega ,F)\, . $ Cela assure la simplicit\'e de $
d(\omega ,F)\, . $ Soit
$ F^\prime \in  {\cal F}(G) $ contenant le sous-groupe $ L^\prime( e)\, . $
L'inclusion $ H \subset  L^\prime( e) $ entra\^\i ne
que $ L^\prime( e) $ est le seul sommet singulier de $ d(\omega ,F^\prime) \,
, $ de sorte que $ F^\prime $ est
dans la classe d'\'equivalence de $ F\, . $ Cela termine la d\'emonstration du
th\'eor\`eme 2.3.3 . \par
\bigskip
\noindent{\bf Corollaire 2.3.4.} {\sl Si \/}$ \omega  \in  \Omega^ 2_{\ell}(
G) $ {\sl poss\`ede un diagramme semi-simple \/}$ d(\omega ,F) $ {\sl qui
est semi-nilpotent, alors \/}$ (M,\omega) $ {\sl poss\`ede un feuilletage
lagrangien
invariant par \/}$ G\, . $ \par
\bigskip
\noindent{\sl Preuve.\/} Le th\'eor\`eme 1.4.1 permet de d\'eformer $
d(\omega ,F) $ en un diagramme simple
$ d(\omega ,F^\prime) \, $ \par
\smallskip
\noindent Un des r\'esultats clefs de ce travail se d\'emontre en conjugant le
th\'eor\`eme
2.3.3 ci-dessus et le th\'eor\`eme 2 de BENSON--GORDON dans $\lbrack$BG
2$\rbrack$. Voici le
r\'esultat clef en question. \par
\bigskip
\noindent{\bf Th\'eor\`eme 2.3.5.} {\sl Soit \/}$ G $ {\sl un groupe de Lie
compl\`etement r\'esoluble contenant un
r\'eseau cocompact \/}$ \Gamma \, . $ {\sl On suppose que la vari\'et\'e
compacte \/}$ \Gamma \backslash G $ {\sl poss\`ede une
m\'etrique Kahlerienne de classe de K\"ahler \/}$ [\alpha] \in
H^2_{DR}(\Gamma \backslash G)\, . $ {\sl Alors toute forme
\/}$ \omega  \in  \Omega^ 2_{\ell}( G) \cap  [\alpha] $ {\sl poss\`ede un
diagramme pond\'er\'e simple \/}$ d(\omega ,F)\  $ \par
\bigskip
\noindent{\sl Esquisse de d\'emonstration.\/} On commence par identifier $
\Omega_{\ell}( G) $ avec $ \Lambda \left(g^\ast \right) = $
l'alg\`ebre ext\'erieure du dual $ g^\ast $ de l'alg\`ebre $ g $ de $ G\, . $
la suite spectrale
de Hattori assure que $ [\alpha]  \cap  \Omega^ 2_{\ell}( G) $ n'est pas vide,
$\lbrack$HA$\rbrack$. Tout $ \omega  \in  [\alpha]  \cap  \Omega^ 2_{\ell}( G)
$
est symplectique. Soit $ {\cal D}G $ le sous-groupe des commutateurs de $ G\,
. $ D'apr\`es
le th\'eor\`eme 2 de $\lbrack$BG 2$\rbrack$ si $ F \in  {\cal F}(G) $ contient
$ {\cal D}G $ alors $ {\cal D}G $ est un sommet
singulier r\'epulsif de poids nul. Le th\'eor\`eme 1.4.1 s'applique \`a cette
situation. Par cons\'equent $ \omega $ poss\`ede un diagramme simple
. \par
\bigskip
\noindent{\sl Remarque 2.3.6.\/} La 2-forme $ \omega  \in  [\alpha]  \cap
\Omega^ 2_{\ell}( G) $ est une forme symplectique. Le
th\'eor\`eme 2.3.3 assure que $ (G,\omega) $ poss\`ede un feuilletage
lagrangien
invariant par les translations \`a gauche. Ce feuilletage se projette en un
feuilletage de $ \Gamma \backslash G\, . $ En g\'en\'eral on n'obtient pas
ainsi un feuilletage
lagrangien pour $ (\Gamma \backslash G,\alpha) \, . $ Si la classe de K\"ahler
$ [\alpha] $ est G-homog\`ene, i.e $ \alpha $
se rel\`eve dans $ G $ en une forme invariante \`a gauche, on aura le
r\'esultat
suivant \par
\bigskip
\noindent{\bf Corollaire 2.3.7} {\sl Les donn\'ees \/}$ G,\Gamma $ {\sl sont
comme dans le th\'eor\`eme 2.3.5. Si \/}$ \Gamma \backslash G $
{\sl poss\`ede une structure K\"alerienne \/}$ (\Gamma \backslash G,\alpha ,J)
$ {\sl dont la forme de K\"ahler \/}$ \alpha $ {\sl est
\/}$ G ${\sl -homog\`ene, alors \/}$ (\Gamma \backslash G,\alpha) $ {\sl
poss\`ede un feuilletage lagrangien }. \par
\bigskip
\noindent La section suivante est domin\'e par des consid\'erations
topologico-g\'eom\'etriques
li\'ees \`a la conjecture des tores plats dont voici l'\'enonc\'e: \par
{\bf Conjecture:} Soit $\Gamma$ un r\'eseau cocompact dans le groupe de lie compl\`etement resoluble $G$. Si la vari\'et\'e $ {\Gamma} \backslash G$ poss\`ede une m\'etrique de kahler alors elle est diff\'eomorphisme au tore plat.\\
$[NB6]$ est consacr\'e partiellement \`a cette conjecture .
\skiplines 2
\noindent{\bf 3. G\'EOM\'ETRIE DES FORMES FERM\'EES.} \par
Cette section est consacr\'ee aux rappels des notions de g\'eom\'etrie affine qui seront
utiles pour une d\'emonstration de la conjecture des tores plats. Certains
\'enonc\'es sont une re\'ecriture adapt\'ee \`a nos objectifs actuels des
r\'esultats
contenus dans les $\lbrack$NB j$\rbrack$, $ j = 1,4,5. $ \par
\noindent Rappelons quelques d\'efinitions utiles dans la suite. \par
\noindent Soit $ (M,\omega) $ une structure symplectique dans la vari\'et\'e $
M\, . $ \par
\bigskip
\noindent{\bf D\'efinition 3.0.1.} {\sl Une structure bilagrangienne dans \/}$
(M,\omega) $ {\sl consiste en deux
feuilletages lagrangiens \/}$ {\cal F} $ {\sl et \/}$ {\cal F}^\prime $ {\sl
qui sont partout transverses }. \par
\bigskip
Soit $ G $ un groupe de Lie compl\`etement r\'esoluble. Soit $ \omega $ une
2-forme
ferm\'ee invariante \`a gauche dans $ G\, . $ Soit $ H $ le sous-groupe de Lie
connexe de $ G $ dont
l'alg\`ebre de Lie engendre le noyau de $ \omega \, . $ Soit $ (M,\omega) $ la
vari\'et\'e
symplectique d\'efinie par $ (G,H,\omega) \, . $ On suppose une fois pour
toutes que $ \omega $
poss\`ede un diagramme simple $ d(\omega ,F)\, , $ dont l'unique sommet
singulier est
not\'e $ S\, . $ Posons $ 2m = \dim M $ et $ \ell  = \dim H\, . $ Les
coordonn\'ees de $ S $ sont la
paire $ \left(G_{m+\ell} ,H_{m+\ell} \right)\, . $ Le sous-groupe $
G_{m+\ell} $ contient les sous-groupes $ H_k $ des sommets
$ S_k= \left(G_k,H_k \right) $ pour $ m + \ell  < k\, . $ Le sous-groupe $
G_{m+\ell} $ engendre dans $ G_k $ le
feuilletage par les classes $ g G_{m+\ell} \, , $ $ g \in  G_k\, , $ $ m +
\ell  < k\, . $ \par
\noindent On va calquer les m\'ethodes mises au point dans les travaux
ant\'erieurs
pour construire dans chaque $ M_k $ une structure bilagrangienne dont un des
feuilletages est engendr\'e par $ G_{m+\ell} $ comme d\'ecrit ci-dessus. (voir
$\lbrack$NB 5$\rbrack$
par exemple). Pour cela on va montrer qu'il y a une relation int\'eressante
entre $ M_k $ et $ M_{k+1}\, . $ \par
\bigskip
{\bf 3.1. }$ {\bf M}_ {\bf k} ${\bf\ comme \lq\lq
M.W}--{\bf r\'eduction\rq\rq\ de }$ {\bf M}_{ {\bf k+1}} {\bf \, .} $ \par
Le but de ce num\'ero 3.1 est de montrer que $ M_k $ peut \^etre construite
comme r\'eduction de MARSDEN--WEINSTEIN de $ M_{k+1}\, . $ Consid\'erons le
graphe $ {\rm gr}(\omega ,F) $
r\'eduit \`a la portion portant les  sommets $ G_k $ et $ G_{k+1} $
$$ \matrix{\displaystyle \leftarrow & \displaystyle H_k & \displaystyle
\leftarrow & \displaystyle H_{k+1} & \displaystyle \leftarrow
\cr\displaystyle  & \displaystyle \downarrow & \displaystyle  & \displaystyle
\downarrow & \displaystyle \cr\displaystyle \rightarrow & \displaystyle G_k &
\displaystyle \rightarrow & \displaystyle G_{k+1} & \displaystyle \rightarrow
\cr\displaystyle  & \displaystyle \downarrow & \displaystyle  & \displaystyle
\downarrow & \displaystyle \cr\displaystyle \leftarrow & \displaystyle M_k &
\displaystyle \leftarrow & \displaystyle M_{k+1} & \displaystyle \leftarrow\cdot \cr} $$
\noindent \overfullrule=0mm \par
\noindent En posant $ {\cal M}_{k+1}= G_k/H_{k+1} $ on obtient le diagramme
commutatif
$$ \matrix{\displaystyle{\cal M}_{k+1} & \displaystyle \doubleup{ i \cr
\hookrightarrow \cr} & \displaystyle M_{k+1} \cr\displaystyle\downarrow \pi &
\displaystyle  & \displaystyle \cr\displaystyle M_k & \displaystyle \cdot &
\displaystyle \cr} \leqno (16) $$
Les formes symplectiques des vari\'et\'es symplectiques $ \left(M_k,\omega_ k
\right) $ sont li\'ees par
$ \pi^ \ast \omega_ k= i^\ast \omega_{ k+1}\, . $ Par ailleurs il est facile
de v\'erifier que $ H_{k+1} $ est un sous-groupe
de Lie distingu\'e de $ H_k $ et que sa codimension ibidem est 1. Soit $ \{
\exp t\zeta ,t\in \bkR\} $
un sous-groupe \`a un param\`etre dans $ H_k $ tel que $ H_k $ est produit
semi-direct
de $ \{ \exp t\zeta ,t\in \bkR\} $ par $ H_{k+1}\, . $ Le champ de vecteurs $
\zeta $ donne lieu \`a un champ
des vecteurs dans $ M_{k+1} $ qui est invariant par $ G_{k+1}\, . $ \par
\noindent Consid\'erons maintenant le champ localement hamiltonien $ \tilde
\zeta $ d\'efini par
$$ \tilde  \zeta( x) = {d \over dt}[\exp t\zeta \cdot x]_{t=0}\, . $$
Naturellement ce champ est partout tangent au feuilletage de $ M_{k+1} $ par
les
hypersurfaces $ \gamma \cdot{\cal M}_{k+1}, \gamma  \in  G_{k+1}\, . $ \par
\noindent Nos hypoth\`eses g\'en\'erales sur la connexit\'e et la simple
connexit\'e des
groupes de Lie consid\'er\'es entra\^\i nent la simple connexit\'e des $
M_k\, . $ La forme
$ i \left(\tilde  \zeta \right)\omega_{ k+1} $ poss\`ede donc une primitive $
\mu_ \zeta \, . $ Il en r\'esulte que l'action de
$ \{ \exp t\zeta ,t\in \bkR\} $ dans $ M_{k+1} $ est hamiltonienne.
L'application moment $ \mu  : M_{k+1}\mapsto  \bkR $
\'etant d\'efinie par
$$ \langle \mu( x),\zeta\rangle  = \mu_ \zeta( x)\, . $$
Il est clair qu'en restriction \`a la sous-vari\'et\'e $ {\cal M}_{k+1} $ la
forme $ i \left(\tilde  \zeta \right)\omega_{ k+1} $
est identiquement nulle. Cela montre que $ {\cal M}_{k+1} $ est une surface de
niveau de
l'application moment. Comme les orbites de $ \{ \exp t\zeta ,t\in \bkR\} $
dans $ {\cal M}_{k+1} $ ne sont
pas autre chose que les feuilles du noyau de $ i^\ast \omega_{ k+1}= \pi^
\ast \omega_ k $ la vari\'et\'e
r\'eduite est la vari\'et\'e $ M_k\ $ \par
\noindent La non \'equivalence de deux sommets adjacents de $ d(\omega ,F) $ a
un contenu
g\'eom\'etri\-que. Ainsi dans la portion de $ d(\omega ,F) $ constitu\'ee des
sommets
hyperboliques d\'ecroissants
$$ \rightarrow  \doubleup{ S^\prime \cr {\rm O} \cr}  ^{\leftarrow}_{ \rightarrow} \doubleup{ S \cr {\rm O} \cr}  \rightarrow $$
la vari\'et\'e symplectique associ\'ee \`a $ S^\prime $ se d\'eduit de celle
associ\'ee au
sommet adjacent $ S $ par r\'eduction de l'espace des phases suivant la proc\'edure de MARSDEN-WEINSTEIN, $ [ {\rm M\ W}]\, . $ Comme cons\'equence de ces relations entre
des sommets
non \'equivalents, nous allons partir d'objets g\'eom\'etriques dans $ M_k $
pour
construire des analogues dans $ M_{k+1}\, . $ \par
\noindent Rappelons que les $ M_k $ de (16) sont d\'eduites du diagramme
pond\'er\'e simple
$ d(\omega ,F)\, . $ Soit $ L_k $ le feuilletage lagrangien de $
\left(M_k,\omega_ k \right) $ engendr\'e par le
sous-groupe $ G_{m+\ell} \, . $ Naturellement $ L_{k+1} $ et $ L_k $ sont
compatibles dans le sens
suivant : \par
\noindent (a) On observe que les feuilles de $ L_{k+1} $ passant par les
points de $ {\cal M}_{k+1} $
sont tout enti\`eres contenues dans $ {\cal M}_{k+1}\, . $ \par
\noindent (b) Les composantes connexes des images inverses par $ \pi  ; {\cal
M}_{k+1}\mapsto  M_k $ des
feuilles de $ L_k $ sont des feuilles de $ L_{k+1}\, . $ \par
\bigskip
\noindent{\bf Th\'eor\`eme }${\bf \lbrack}${\bf NB 5}${\bf \rbrack}${\bf .}
{\sl On suppose que \/}$ M_k $ {\sl porte un deuxi\`eme feuilletage
lagrangien \/}$ {\cal N}_k $ {\sl transverse \`a \/}$ L_k\, . $ {\sl Alors
\/}$ M_{k+1} $ {\sl porte aussi un deuxi\`eme
feuilletage lagrangien \/}$ {\cal N}_{k+1} $ {\sl transverse \`a \/}$ L_{k+1}
$ {\sl et compatible avec \/}$ {\cal N}_k $ {\sl dans
le sens suivant : les composantes des traces dans \/}$ {\cal M}_{k+1} $ {\sl
des feuilles de \/}$ {\cal N}_{k+1} $
{\sl sont des rev\^etements des feuilles de \/}$ {\cal N}_k\  $ \par
\bigskip
\noindent{\sl Esquisse de la d\'emonstration.\/} On va se restreindre \`a
d\'ecrire les \'etapes de
la d\'emonstration de fa\c con \`a ce que le lecteur puisse \'eventuellement
v\'erifier les d\'etails sans difficult\'es violentes. \par
\bigskip
\noindent{\bf 1\`ere \'etape.} On choisit un sous-groupe \`a un param\`etre $
\{ \exp t\xi ,t\in \bkR\}  \subset  G_{k+1} $
tel que $ G_{k+1}= G_k\times  \{ \exp t\xi ,t\in \bkR\} $ le champ de vecteurs
hamiltonien $ \tilde  \xi $ est sans
point fixe dans $ M_{k+1}\, . $ Le syst\`eme diff\'erentiel $ i \left(\tilde
\xi \right)\omega_{ k+1}= 0 $ est
compl\`etement int\'egrable. On d\'esigne par $ {\cal G}_k $ le feuilletage de
$ M_{k+1} $ par les
$ \gamma \cdot{\cal M}_{k+1}, $ $ \gamma \in G_{k+1}\, . $ L'intersection de $
{\cal G}_k $ avec $ \gq\gq i \left(\tilde  \xi \right)\omega_{ k+1}= 0\dq\dq $
est un feuilletage
de codimension 2 dans $ M_{k+1} $ subordonn\'e \`a $ {\cal G}_k\  $
\par
\bigskip
\noindent{\bf 2\`eme \'etape.} On consid\`ere le feuilletage de $ {\cal
M}_{k+1} $ par les feuilles de
$ {\cal G}_k\cap  \left[i \left(\tilde  \xi \right)=0 \right]\, . $ On note $
{\cal D} $ la distribution tangente ce feuilletage de $ {\cal M}_{k+1}\, ; $
$ {\cal D} $ est horizontal pour la projection $ \pi  : {\cal M}_{k+1}\mapsto
M_k $ dans le sens suivant : $ d\pi $
envoie bijectivement $ {\cal D}_y $ sur $ T_{\pi( y)}M_k\, . $ Il existe dans
$ {\cal D} $ une
sous-distribution $ {\cal D}_k $ qui est projet\'ee par $ \pi $ sur $ {\cal
N}_k\, . $ Ainsi les feuilles de
$ {\cal D}_k $ sont des rev\^etements des feuilles de $ {\cal N}_k\, . $
Evidemment $ {\cal D}_k $ est transverse \`a
$ L_{k+1} $ dans $ {\cal M}_{k+1}\, . $ \par
\bigskip
\noindent{\bf 3\`eme \'etape.} On observe que le champ de vecteurs $ \tilde
\xi $ est sans z\'ero dans $ M_{k+1} $
il d\'efinit donc un feuilletage par des courbes r\'eguli\`eres. De plus
chaque
courbe int\'egrale de $ \tilde  \xi $ est une transversale pour le feuilletage
$ {\cal G}_k $ de $ M_{k+1}\, . $
On peut donc construire un \lq\lq feuilletage r\'egl\'e\rq\rq\ de $ M_{k+1} $
par glissement des
feuilles de $ {\cal D}_k $ le long des courbes int\'egrales de $ \tilde  \xi
\, . $ Plus pr\'ecis\'ement si
$ N(x) $ est une feuille de $ {\cal D}_k $ passant par $ x \in  {\cal M}_{k+1}
$ alors on d\'esigne par $ {\cal N}_{k+1}(x) $
la sous-vari\'et\'e de $ M_{k+1} $ d\'efinie par
$ {\cal N}_{k+1}(x) = \{ \exp t\xi \cdot( N(x)),t\in \bkR\}  = \doublelow{
\cup \cr t\in \bkR \cr}( \exp t\xi)  (N(x))\, . $ C'est une
sous-vari\'et\'e connexe de dimension $ = {1 \over 2} \dim M_{k+1}\, , $ de
plus $ {\cal N}_{k+1}(x) $ est
lagrangien. \par
\bigskip
\noindent{\bf 4\`eme \'etape.} Il est facile de voir que si $ N(x)\cap
N(x^\prime) $ est vide, alors il en
sera de m\^eme de $ {\cal N}_{k+1}(x) \cap  {\cal N}_{k+1}(x^\prime) $ et que
$ M_{k+1} $ est r\'eunion des $ {\cal N}_{k+1}(x)\, . $ \par
\noindent Pour terminer on observera que le feuilletage $ {\cal N}_{k+1} $
ainsi construit est
partout transverse \`a $ L_{k+1} $ et poss\`ede la propri\'et\'e de
compatibilit\'e avec $ {\cal N}_k $
qui est attendue. Ici se termine la d\'emonstration du th\'eor\`eme
$\lbrack$NB 5$\rbrack \, $ \par
\smallskip
Maintenant il est facile de montrer que celle des vari\'et\'es $
\left(M_k,\omega_ k \right) $
qui est de dimension 2 poss\`ede une structure bilagrangienne dont l'un des
deux feuilletages est engendr\'e par $ G_{m+\ell} \, . $ Le th\'eor\`eme
$\lbrack$NB 5$\rbrack$ ci-dessus
assure qu'il est de m\^eme dans chacune des $ \left(M_k,\omega_ k \right)\ $ \par
\bigskip
{\bf 3.2. Structures bilagrangiennes affines }${\bf
\lbrack}${\bf NB 1}${\bf \rbrack}${\bf , }${\bf \lbrack}${\bf NB 5}${\bf
\rbrack}${\bf , }${\bf \lbrack}${\bf HH}${\bf \rbrack}${\bf .} \par
Lorsque la vari\'et\'e symplectique $ (M,\omega) $ poss\`ede deux feuilletages
lagrangiens transverses $ {\cal D}_1 $ et $ {\cal D}_2 $ il existe une unique
connexion lin\'eaire $ D $
satisfaisant les conditions suivantes \par
\smallskip
\noindent (i) Le tenseur de torsion de $ D $ est nul ; \par
\smallskip
\noindent (ii) $ D\omega  = 0\, ; $ \par
\smallskip
\noindent (iii) Pour tout champ de vecteurs tangent $ X, $ $ D_X \left({\cal
D}_j \right) \subset  {\cal D}_j\, , $ $ j = 1,2\, . $ \par
\noindent Voici l'expression de $ D\, . $ Pour deux champs de vecteurs $ X $
et $ Y $ on d\'efinit $ D^o_XY $
en posant
$$ i \left(D^o_XY \right)\omega  = L_Xi(Y)\omega \, . $$
Posons maintenant $ X = \left(X_1,X_2 \right)\in{\cal D}_1\times  {\cal
D}_2\, ; $ l'expression de $ D_XY $ est :
$$ D_{ \left(X_1,X_2 \right)} \left(Y_1,Y_2 \right) = \left(D^o_{X_1}Y_1+
\left[X_2,Y_1 \right]_1,D^o_{X_2}Y_2+ \left[X_1,Y_2 \right]_2 \right)\, .
\leqno (17) $$
$ ([X,Y]_j $ signifie la $ {\cal D}_j $-composante de $ [X,Y]\, , $ $ j =
1,2). $ \par
\noindent En g\'en\'eral le tenseur de courbure de $ D $ est non nul
. \par
\bigskip
\noindent{\bf Th\'eor\`eme }${\bf \lbrack}${\bf HH}${\bf \rbrack}${\bf .}
{\sl Les conditions suivantes sont \'equivalentes :\/} \par
\noindent{\sl (i) Le tenseur de courbure de \/}$ D $ {\sl est nul.\/} \par
\noindent{\sl (ii) Dans un voisinage de tout point il existe des coordonn\'ees
de Darboux
\/}$ \left(q_1,...,q_m,p_1,...,p_m \right) $ {\sl telle que les \/}$ q_j $
{\sl sont des int\'egrales premi\`eres de \/}$ {\cal D}_1 $
{\sl et les \/}$ p_j $ {\sl sont des int\'egrales premi\`eres de \/}$ {\cal
D}_2\  $ \par
\bigskip
\noindent Dans $\lbrack$NB 4$\rbrack$ on a prouv\'e une version plus forte du
th\'eor\`eme de HESS
ci-dessus. Avant d'\^etre plus pr\'ecis retournons quelques instants \`a la
situation du diagramme (16) :
$$ \matrix{\displaystyle H_k & \displaystyle \leftarrow & \displaystyle
H_{k+1} \cr\displaystyle \downarrow & \displaystyle  & \displaystyle
\downarrow \cr\displaystyle G_k & \displaystyle \rightarrow & \displaystyle
G_{k+1} \cr\displaystyle \downarrow & \displaystyle  & \displaystyle
\downarrow \cr\displaystyle{\cal M}_{k+1} & \displaystyle \doubleup{ i \cr
\hookrightarrow \cr} & \displaystyle M_{k-1} \cr\displaystyle\downarrow \pi &
\displaystyle  & \displaystyle \cr\displaystyle M_k & \displaystyle \cdot &
\displaystyle \cr} $$
On \'ecrit comme pr\'ec\'edemment $ H_k $ sous la forme de produit semi-direct
$$ H_k= H_{k+1}\times  \{ \exp t\xi ,t\in \bkR\} \, . $$
Soit $ x_o\in  M_{k+1} $ le point image de $ H_{k+1} $ par la projection
canonique de $ G_{k+1} $
sur $ M_{k+1}\, . $ Posons $ \zeta_ o= {d \over dt} \exp t\zeta \cdot
\left(x_o \right)_{\mid t=0}\, . $ On d\'efinit pour tout $ x \in  M_{k+1} $
$$ \zeta( x) = d\gamma_{ x_o} \left(\zeta_ o \right) $$
o\`u $ \gamma $ est un \'el\'ement de $ G_{k+1} $ tel que $ x = \gamma
\left(x_o \right)\, . $ Le champ des vecteurs
tangents $ x \mapsto  \zeta( x) $ est bien d\'efini, sans point fixe et
invariant par $ G_{k+1}\, . $
Le champ $ \zeta $ est tangent \`a $ {\cal M}_{k+1} $ et engendre le noyau de
$ \pi^ \ast \omega_ k= i^\ast \omega_{ k+1}\, . $ \par
\noindent Posons aussi
$$ G_{k+1}= G_k\times  \{ \exp t\xi ,t\in \bkR\} $$
Le champ hamiltonien $ x \mapsto  \tilde  \xi( x) $ d\'eduit de l'action de $
\{ \exp t\xi ,t\in \bkR\} $ dans $ M_{k+1} $
est sans z\'ero. La fonction
$$ x \mapsto  \omega_ x \left(\tilde  \xi ,\zeta \right) $$
est partout non nulle. \par
\noindent On d\'efinit le champ de vecteur colin\'eaire \`a $ \zeta ,\hat
\zeta \, , $ en posant en $ x \in  M_{k+1} $
$$ \hat  \zeta( x) = {1 \over \omega_ x \left(\tilde  \xi ,\zeta \right)}
\zeta( x)\, . $$
\par
\medskip
\noindent{\bf Lemme (}${\bf \lbrack}${\bf NB j}${\bf \rbrack}${\bf , j }${\bf
=}${\bf\ 4,5).} {\sl Le champ \/}$ \hat  \zeta $ {\sl est hamiltonien.\/} \par
\smallskip
\noindent Le lemme ci-dessus intervient dans la preuve du th\'eor\`eme
suivant. \par
\bigskip
\noindent{\bf Th\'eor\`eme }${\bf \lbrack}${\bf NB 4}${\bf \rbrack}${\bf .}
{\sl On se place dans les conditions de (16) et on admet que
\/}$ \left(M_k,\omega_ k \right) $ {\sl et \/}$ \left(M_{k+1},\omega_{ k+1}
\right) $ {\sl poss\`edent des structures bilagrangiennes \/}$ \left(L_k,
{\cal N}_k \right) $ {\sl et
\/}$ \left(L_{k+1},{\cal N}_{k+1} \right) $ {\sl qui sont compatibles et dont
\/}$ L_k $ {\sl et \/}$ L_{k+1} $ {\sl sont engendr\'es par
\/}$ G_{m+\ell} \, . $ \par
\smallskip
\noindent{\sl (i) On suppose que dans un voisinage de chaque \/}$ y \in  M_k $
{\sl il existe un syst\`eme
de coordonn\'ees de Darboux \/}$ \left(q_1,...,q_m,p_1,...,p_m \right) $ {\sl
v\'erifiant les conditions
de HESS c'est \`a dire : les feuilles de \/}$ L_k $ {\sl sont d\'efinies par
\/}$ q_j= {\rm  cte} $ {\sl et
celles de \/}$ {\cal N}_k $ {\sl par \/}$ p_j=\ {\rm cte} \, . $ \par
\smallskip
\noindent{\sl (ii) On suppose en outre que si \/}$ \left(q^{\prime}_
1,...,q^{\prime}_ m,p^{\prime}_ 1,...,p^{\prime}_ m \right) $ {\sl est un
autre syst\`eme
des coordonn\'ees de Darboux qui v\'erifie (i), la matrice jacobienne \/}$ {D
\left(q^{\prime}_ 1,...,q^{\prime}_ m \right) \over D \left(q_1,...,q_m
\right)} $
{\sl est constante et appartient au groupe orthogonal \/}$ O(m)\, . $ {\sl
Alors les m\^emes
conclusions (i) et (ii) valent dans un voisinage de tout point \/}$ x \in
M_{k+1} $
pour la paire $ \left(L_{k+1},{\cal N}_{k+1} \right)\, . $ \par
\bigskip
\noindent{\sl Esquisse de preuve.\/} \par
\noindent 1\nobreak$^\circ$) Il suffit de prouver le th\'eor\`eme pour les
points $ x $ de $ {\cal M}_{k+1}\, . $ \par
\bigskip
\noindent{\bf 1\`ere \'etape.} On fixe un voisinage $ U $ de $ x \in  {\cal
M}_{k+1} $ dont l'image par $ \pi  : {\cal M}_{k+1}\mapsto  M_k $
est un ouvert domaine de coordonn\'ees de Darboux $
\left(q_1,...,q_m,p_1,...,p_m \right) $
v\'erifiant (i) et (ii). On remonte les fonctions $ q_j $ et $ p_j $ dans
l'ouvert $ U\, . $ \par
\noindent On obtient ainsi dans $ U $ $ 2m $ fonctions $ \tilde
q_1,...,\tilde  q_m,\tilde  p_1,...,\tilde  p_m\, . $ Pour $ \varepsilon  > 0
$
assez petit $ \doublelow{ \cup \cr\vert t\vert  < \varepsilon \cr} \exp(t\xi)
\cdot U $ est un voisinage de $ x $ dans $ M_{k+1}\, . $ \par
\noindent On prolonge les fonctions $ \tilde  q_j $ et $ \tilde  p_j $ dans $
U_\varepsilon = \doublelow{ \cup \cr\vert t\vert  < \varepsilon \cr}
\exp(t\xi) \cdot \cup $ en posant
$$ \tilde {q}_{j}((\exp t\xi).y) = \tilde {q}_{j}(y) \quad pour \quad y \in U , $$
$$ \tilde {p}_{j}((\exp t\xi).y) = \tilde {p}_{j}(y) . $$

Quitte \`a \lq\lq diminuer\rq\rq\ $ U_\xi \, , $ on peut y trouver deux
fonctions diff\'erentiables
$ q_{m+1} $ et $ p_{m+1} $ v\'erifiant
$$ i \left(\hat  \zeta \right)\omega_{ k+1}= dq_{m+1}\ {\rm et} \ i
\left(\tilde  \xi \right)\omega_{ k+1}= dp_{m+1} $$
Puisque $ \omega \left(\tilde  \xi ,\hat  \zeta \right) = 1\, , $ le syst\`eme
des coordonn\'ees locales
$$ \left(\tilde  q_1,...,\tilde  q_m,q_{m+1},\tilde  p_1,...,\tilde
p_m,p_{m+1} \right) $$
v\'erifie (i). Naturellement les syst\`emes des coordonn\'ees dans $ M_{k+1} $
obtenus
comme ci-dessus v\'erifient les conditions (i) et (ii) . \par
\smallskip
\noindent Il est facile de voir que celles de $ \left(M_k,\omega_ k \right) $
qui est de dimension 2 poss\`ede
un atlas de Darboux v\'erifiant les conclusions (i) et (ii) du th\'eor\`eme
ci-dessus. Il suffit pour cela de prendre des applications moment des
actions de $ \hat  \zeta $ et de $ \tilde  \xi $ respectivement. On d\'eduit
du th\'eor\`eme $\lbrack$NB 4$\rbrack$
ci-dessus que la vari\'et\'e symplectique $ (M,\omega) $ d'o\`u se d\'eduisent
les $ \left(M_k,\omega_ k \right) $
par r\'eduction symplectiques successives poss\`ede un atlas de Darboux
v\'erifiant les conditions (i), (ii) pour la structure bilagrangienne $
(L,{\cal N})\  $ \par
\medskip
\noindent{\sl Remarques 3.2.1.\/} (a) En vertu du th\'eor\`eme de HESS la
lcondition (i)
\'equivaut \`a la nullit\'e du tenseur de courbure de la connexion d\'efinie
dans
l'expression (17). \par
\smallskip
\noindent (b) La condition (ii) de changement de coordonn\'ees locales montre
que les
fonctions \`a valeurs complexes
$$ z_j= q_j+ \sqrt{ -1} p_j $$
satisfont l'\'equation de Cauchy--Riemann. Elles d\'efinissent par
cons\'equent
une structure de vari\'et\'e analytique complexe dans $ M_k\, . $ \par
\smallskip
\noindent (c) On d\'eduit \'egalement de la conditions (ii) que les formes
hermitiennes
locales
$$ h = \sum^{ }_ jdz_jd\bar  z_j $$
satisfont la condition de recollement. Autrement dit, sur l'intersection
des domaines de d\'efnition des $ z_j $ et des $ z^{\prime}_ j $ on a
$$ \sum^{ }_ jdz_jd\bar  z_j= \sum^{ }_{ } dz^{\prime}_ jd\bar  z^{\prime}_
j\, . $$
Par cons\'equent les $ M_k $ sont des vari\'et\'es kahleriennes. De plus les
structures kahleriennes de $ M_k $ et de $ M_{k+1} $ sont compatibles dans le
sens
suivant : si $ \left(z_1...z_m \right) $ sont des coordonn\'ees locales
holomorphes de $ M_k\, , $
alors leurs images inverses $ \pi^ \ast z_j $ sont des restrictions \`a $
{\cal M}_{k+1} $ des
fonctions analytiques complexes locales d\'efinies dans $ M_{k+1}\, . $ \par
\smallskip
\noindent (d) On montre en outre que la connexion lin\'eaire $ D $ d\'efinie
dans (17) est
la connexion de Levi--Civita de $ M_k $ munie de la structure k\"ahlerienne
ci-dessus . \par
\noindent Nous allons donner un condens\'e des r\'esultats mis en \'evidence
ci-dessus
sous la forme suivante. \par
\bigskip
\noindent{\bf Th\'eor\`eme 3.2.2. }{\sl Soit \/}$ \omega  \in  \Omega^
2_{\ell}( G) $ ;$ H $ {\sl est le sous-groupe
de Lie connexe de \/}$ G $ {\sl dont l'alg\`ebre de Lie engendre le noyau de
\/}$ \omega \, ; $ $ (M,\omega) $
{\sl est la vari\'et\'e symplectique \/}$ G ${\sl -homog\`ene d\'etermin\'ee
par le triplet \/}$ (G,H,\omega) \, . $
{\sl On suppose que la forme \/}$ \omega $ {\sl poss\`ede un diagramme simple
\/}$ d(\omega ,F)\, . $ {\sl Alors la
vari\'et\'e symplectique \/}$ (M,\omega) $ {\sl poss\`ede une structure
kahlerienne affinement
plate \/}$ (M,\omega ,J) $ {\sl et un feuilletage lagrangien \/}$ L $ {\sl
dont les feuilles sont des
sous-vari\'et\'es totalement r\'eelles de la vari\'et\'e analytique complexe
\/}$ (M,J)\,  $ \par
\bigskip
\noindent{\sl Remarques 3.2.3.\/} (a) La structure Kahlerienne $ (M,\omega ,J)
$ est celle d\'etermin\'ee
par une structure bilagrangienne $ (L,{\cal N}) $ et la structure affine de $
M $ est
celle d\'efinie par la connexion localement $ D $ du couple $ (L,{\cal N}) $
suivant
l'expression (17). \par
\smallskip
\noindent (b) Le th\'eor\`eme 3.2.2 dit que $ D $ est la connexion de
Levi-Civita de la
m\'etrique riemannienne $ g_J $ d\'efinie par $ g_J(u,v) = \omega( Ju,v)\ $ \par
\smallskip
\noindent La connexion lin\'eaire $ D $ pr\'eserve $ L $ et $ {\cal N} $ d'une
part et d'autre part on a
$ DJ = 0 $ et $ D\omega  = 0\, . $ Ainsi du point de la g\'eom\'etrie
riemannienne de $ \left(M,g_J \right) $
les feuilles de $ L $ (resp. de $ {\cal N}) $ sont des sous-vari\'et\'es
totalement
g\'eod\'esiques. \par
\smallskip
\noindent (c) Du point de vue de la g\'eom\'etrie affine du couple $ (M,D) $
les feuilles de
$ L $ (resp. de $ {\cal N}) $ sont les sous-vari\'et\'es affines .
\par
\smallskip
\noindent (d) On observe que des structures mises en \'evidence dans $ M $ les
seules
visiblement $ G $-homog\`enes sont la vari\'et\'e symplectique $ (M,\omega) $
et la vari\'et\'e
feuillet\'ee $ (M,L)\, . $ \par
\smallskip
\noindent (e) Enfin puisque $ D\omega  = 0 $ l'holonomie lin\'eaire de $ D $
est un sous-groupe
discret de $ SL(m,\bkR) \, , $ $ m = \dim M\, . $ On sait que dans le cas des
vari\'etr\'es
compactes la compl\'etude de $ D $ est conjectur\'ee \^etre intimement li\'ee
\`a cette
inclusion du groupe d'holonomie lin\'eaire dans $ SL(m,\bkR) , $
$\lbrack$FGH$\rbrack$, c'est la
conjecture de MARKUS.\par
\skiplines 4
\noindent \overfullrule=0mm \par
\centerline{{\bf R\'EF\'ERENCES}}
\bigskip
\noindent${\bf \lbrack}${\bf ABCKT}${\bf \rbrack}${\bf\ J. AMOROS, N. BURGER,
K. CORLETTE, D. KOSCHICK and D. TOLEDO.
}{\sl Fundamental groups of compact Kahler manifolds.\/}{\bf\ }Amer Math. Soc.
Math. survey
and Monographs vol. 44. \par
\smallskip
\noindent${\bf \lbrack}${\bf AU 1}${\bf \rbrack}${\bf\ L. AUSLANDER. }{\sl
Bieberbach's groups on space groups and discrete
uniform subgroups of Lie groups }$\Irm${\sl . \/}Ann. Math. 71 (1960)
579--590. \par
\smallskip
\noindent${\bf \lbrack}${\bf AU 2}${\bf \rbrack}${\bf\ L. AUSLANDER. }{\sl An
exposition of the structures of solvmanifolds. \/}Bull.
of Amer. Math. Soc. 79 (1973) 227--261. \par
\smallskip
\noindent${\bf \lbrack}${\bf BG 1}${\bf \rbrack}${\bf\ C. BENSON and C.
GORDON. }Kahler and symplectic structures on
nilmanifolds topology 27 (1988) 513--518. \par
\smallskip
\noindent${\bf \lbrack}${\bf BG 2}${\bf \rbrack}$ {\bf C. BENSON and C.
GORDON.} {\sl Kahler structures on compact
solvmanifolds. \/}Proc. of the Amer. Math. Soc. 108 (1990) 971--990. \par
\smallskip
\noindent${\bf \lbrack}${\bf BE}${\bf \rbrack}${\bf\ A. BESSE.}{\sl\ Enstein
manifolds.\/} Ergebnisse der Math. und ihrer
Grenzgebiete 3. Folge. Band 10. \par
\smallskip
\noindent${\bf \lbrack}${\bf B0}${\bf \rbrack}${\bf\ N. BOURBAKI.} {\sl
Groupes et Alg\`ebres de Lie.\/} vol. 7. Hermann Paris. \par
\smallskip
\noindent${\bf \lbrack}${\bf BW}${\bf \rbrack}${\bf\ F. BRUHAT and H.
WHITNEY.} Comment. Math. Helv. 33 (1959) 132--160. \par
\smallskip
\noindent${\bf \lbrack}${\bf CFG}${\bf \rbrack}${\bf\ L. CORDERO, M. FERNANDEZ
and A. GRAY. }{\sl Symplelctic manifolds with no
Kahler structure.\/} Topology 25 (1986) 375--380. \par
\smallskip
\noindent${\bf \lbrack}${\bf DGMS}${\bf \rbrack}${\bf\ P. DELIGNE, P.
GRIFFITHS, J. MORGAN and D. SULLIVAN. }{\sl Real homotopy
theory of Kahler manifolds.\/} Inv. Math. 29 (1975) 245--274. \par
\smallskip
\noindent${\bf \lbrack}${\bf DN}${\bf \rbrack}${\bf\ J. DORFNEISTER and K.
NAKAJIMA.} {\sl The fundamental conjecture Acta Math.
160 (1988) 13--70.\/} \par
\smallskip
\noindent${\bf \lbrack}${\bf FGH}${\bf \rbrack}${\bf\ D. FRIED, W. GOLDMAN and
M.W. HIRSCH. }{\sl Affine manifolds with
nilpotent holonomy group.\/} Comment Math. Helv. 56 (1981) 487--523. \par
\smallskip
\noindent${\bf \lbrack}${\bf G}${\bf \rbrack}${\bf\ M. GROMOV.} {\sl Sur le
groupe fondamental d'une vari\'et\'e kahl\'erienne.\/} CR.
Acad. Sci. Paris 308 (1989) 67--70. \par
\smallskip
\noindent${\bf \lbrack}${\bf GPV 1}${\bf \rbrack}${\bf\ S.G. GINDIKIN, I.I.
PIATECCKII, SAPIRO and E.B. VINBERG}.{\bf\
}{\sl Homogeneous K\"ahler manifolds in \lq\lq Geometry of homogeneous bounded
domains\rq\rq .\/}
CIME 1967. \par
\smallskip
\noindent${\bf \lbrack}${\bf GPV 2}${\bf \rbrack}${\bf\ S.G. GINDIKIN, I.I.
PIATECCKII, SAPIRO and E.B. VINBERG. }{\sl On the
classification and canonical realization of complex homogeneous bounded
domains.\/} Trans. Moscow Math. Soc. 13 (1964) 340--403. \par
\smallskip
\noindent${\bf \lbrack}${\bf GH}${\bf \rbrack}${\bf\ W. GOLDMAN and M.W.
HIRSCH. }{\sl Affine manifolds and Orbits of Algebraic
groups.\/} Trans. Amer. Math. Soc. 295 (1986) 175--197. \par
\smallskip
\noindent${\bf \lbrack}${\bf GM}${\bf \rbrack}${\bf\ W. GOLDMAN and J.J.
NILLSON.} {\sl The deformation theory of
representations of fundamental groups of compact K\"ahler manifolds.\/} Publ.
Math. IHES 67 (1988) 43--96. \par
\smallskip
\noindent${\bf \lbrack}${\bf GS}${\bf \rbrack}${\bf\ V. GUILLEMIN and M.
STENZEL. }{\sl Grauert tube and the homogeneous
Monge-Ampere equations.\/} J. Diff. Geom. 34 (1991) 561--570. \par
\smallskip
\noindent${\bf \lbrack}${\bf HA}${\bf \rbrack}${\bf\ A. HATTORI. }{\sl
Spectral sequence in the de RHAM cohomology of fibre
bundles.\/} Fac. sei. Univ. Tokyo 8 (1960) 289--331. \par
\smallskip
\noindent${\bf \lbrack}${\bf HH}${\bf \rbrack}${\bf\ A.HESS. }{\sl Connection
on symplectic manifolds and Geometric
quantization.\/} Lecture Notes in Math. \par
\smallskip
\noindent${\bf \lbrack}${\bf HW}${\bf \rbrack}${\bf\ W.Y. HISIANG. }{\sl
Cohomology theory of Topological Transf. groups Erg.
der Math. une chrer Gr. band 85\/} Springer--Verlag. \par
\smallskip
\noindent${\bf \lbrack}${\bf KL 1}${\bf \rbrack}${\bf\ J.L. KOSZUL. }{\sl
Homologie et Cohomologir des Alg\`ebres de Lie.\/} Bull.
Soc. Math. France 48 (1950) 65--127. \par
\smallskip
\noindent${\bf \lbrack}${\bf KL 2}${\bf \rbrack}${\bf\ J.L. KOSZUL. }{\sl
D\'eformations des vari\'et\'es localement plates.\/} Ann.
Inst. Fourier 18 (1968). \par
\smallskip
\noindent${\bf \lbrack}${\bf L}${\bf \rbrack}${\bf\ P. LIBERMANN. }{\sl
Probl\`eme d'\'equivalence en g\'eom\'etrie symplectique.\/}
Ast\'erisque 107 (1983) 43--68. \par
\smallskip
\noindent${\bf \lbrack}${\bf MA}${\bf \rbrack}${\bf\ Y. MATSUSHIMA.} {\sl On
discrete subgroups and Homogeneous spaces of
nilpotent Lie groups.\/} Nagoya Math. Jour. 2 (1951) 95--110. \par
\smallskip
\noindent${\bf \lbrack}${\bf Mc D 1}${\bf \rbrack}${\bf\ Dusa Mc DUFF.}{\sl\
The moment map for circle action on symplectic
manifolds JGP 5 (1988) 149--160.\/} \par
\smallskip
\noindent${\bf \lbrack}${\bf Mc D 2}${\bf \rbrack}${\bf\ Dusa Mc DUFF. }{\sl
Example of simply connected symplectic manifolds
non k\"ahlerian.\/} Jour. Diff. Geom. 70 (1984) 267--277. \par
\smallskip
\noindent${\bf \lbrack}${\bf MI}${\bf \rbrack}${\bf\ J. MILNOR.} {\sl
Fundamental groups of complete affinely flat manifolds.\/}
Adv. in Math. 25 (1977) 178--187. \par
\smallskip
\noindent${\bf \lbrack}${\bf MO 1}${\bf \rbrack}${\bf\ G.D. MOSTON.} {\sl
Cohomology of topological groups and solomanifolds.\/}
Ann. of Math. 73 (1961) 20--48. \par
\smallskip
\noindent${\bf \lbrack}${\bf MO 2}${\bf \rbrack}${\bf\ G.D. MOSTON. }{\sl
Fundamental groups of homogeneous spaces.\/} Ann. of
Math. 66 (1957) 249--255. \par
\smallskip
\noindent${\bf \lbrack}${\bf MW}${\bf \rbrack}${\bf\ J. MARSDEN and WEINSTEIN.
}{\sl Reduction of symplectic manifolds with
symmetry.\/} Reports on Math. Phys. 5 (1974) 121--130. \par
\smallskip
\noindent${\bf \lbrack}${\bf NB 1}${\bf \rbrack}${\bf\ NGUIFFO B. BOYOM.
}{\sl Vari\'et\'es symplectiques affines.\/} Manuscripta Math.
64 (1989) 1--33. \par
\smallskip
\noindent${\bf \lbrack}${\bf NB 2}${\bf \rbrack}${\bf\ NGUIFFO B. BOYOM.
}{\sl Structures affines homotpes \`a z\'ero.\/} Jour. Diff.
Geom. 33 (1990) 859--911. \par
\smallskip
\noindent${\bf \lbrack}${\bf NB 3}${\bf \rbrack}${\bf\ NGUIFFO B. BOYOM.
}{\sl Structures affines isotropes des groupes de Lie.\/}
Annalli Scuo. Norm. Sup. Pisa. \par
\smallskip
\noindent${\bf \lbrack}${\bf NB 4}${\bf \rbrack}${\bf\ NGUIFFO B. BOYOM.}
{\sl M\'etriques k\"ahl\'eriennes affinement plates.\/} Proc.
London Math. Soc. 3 (1993) 358--380. \par
\smallskip
\noindent${\bf \lbrack}${\bf NB 5}${\bf \rbrack}${\bf\ NGUIFFO B. BOYOM.
}{\sl Structures localement plates dans certaines
vari\'et\'es symplectiques.\/} Math Scand. 76 (1995) 61--84. \par
\smallskip
\noindent${\bf \lbrack}${\bf NB 6}${\bf \rbrack}${\bf\ NGUIFFO B. BOYOM.
}{\sl Conjecture des tores plats \/} (Soumis ). \par
\smallskip
\noindent${\bf \lbrack}${\bf NK}${\bf \rbrack}${\bf\ K. NOMIZU.} {\sl On the
cohomology of compact homogeneous spaces of
nilpotent Lie groups.\/} Ann. of Math. 59 (1954) 531--538. \par
\smallskip
\noindent${\bf \lbrack}${\bf ON}${\bf \rbrack}${\bf\ A.L. ONISCIK.} {\sl
Inclusion relations among transitive compact groups.\/}
Amer. Math. Soc. Translations S.2 50 (1966) 5--58. \par
\smallskip
\noindent${\bf \lbrack}${\bf RA}${\bf \rbrack}${\bf\ S. RAGHUNATHAN. }{\sl
Discrete subgroups of Lie groups.\/} Ergebnisse der
Math. 68 Springer Verlag 1972. \par
\smallskip
\noindent${\bf \lbrack}${\bf RC}${\bf \rbrack}${\bf\ C. RAMANUJAM.} {\sl A
topological characterisation of the affine plane as
an Algebraic Variety.\/} Ann. of Math. 94 (1971) 69--88. \par
\smallskip
\noindent${\bf \lbrack}${\bf R.W}${\bf \rbrack}${\bf\ J. ROELS and A.
WEINSTEIN. }{\sl On functions whose Poisson brackets are
constant.\/} J. Math. Phys. 12. 1482--1486. \par
\smallskip
\noindent${\bf \lbrack}${\bf SA}${\bf \rbrack}${\bf\ M. SATO. }{\sl Sur
certains groupes de Lie r\'esolubles. \/}Sci. pap. Univ. Tokyo
7 (1957) $\Irm$--$\IIrm$ ; 157--168. \par
\smallskip
\noindent${\bf \lbrack}${\bf SE}${\bf \rbrack}${\bf\ J.P. SERRE}. {\sl
GAGA.\/} Ann. Inst. Fourier 6 (1956) 1--42. \par
\smallskip
\noindent${\bf \lbrack}${\bf TH}${\bf \rbrack}${\bf\ W. THURSTON. }{\sl Some
simple of symplectic manifolds.\/} Proc. Amer. Math.
Soc. 55 (1976) 456--478. \par
\smallskip
\noindent${\bf \lbrack}${\bf WL}${\bf \rbrack}${\bf\ A. WEILL.} {\sl
Introduction \`a l'\'etude des vari\'et\'es k\"ahl\'eriennes.\/} Hermann
Paris 1958. \par
\smallskip
\noindent${\bf \lbrack}${\bf WN}${\bf \rbrack}${\bf\ A. WEINSTEIN.} {\sl
Symplectic manifolds and their lagrangian submanifolds.\/}
Adv. In Math. 6 (1971) 329--346. \par
\smallskip
\noindent${\bf \lbrack}${\bf WO}${\bf \rbrack}${\bf\ N.M.J. WOODHOUSE.} {\sl
Geometric quantization.\/} Oxford Math. Monographs.
Larendon Press Oxford. \par
\skiplines 4

\end{document}